%% file: strol.tex
%

\magnification=\magstep1
\def\forces{\parallel\!\!\! -}


\def\hexnumber#1{\ifcase#1 0\or1\or2\or3\or4\or5\or6\or7\or8\or9\or
	A\or B\or C\or D\or E\or F\fi }

\font\teneuf=eufm10
\font\seveneuf=eufm7
\font\fiveeuf=eufm5
\newfam\euffam
\textfont\euffam=\teneuf
\scriptfont\euffam=\seveneuf
\scriptscriptfont\euffam=\fiveeuf
\def\frak{\fam\euffam \teneuf}

\font\tenmsx=msam10
\font\sevenmsx=msam7
\font\fivemsx=msam5
\font\tenmsy=msbm10
\font\sevenmsy=msbm7
\font\fivemsy=msbm5
\newfam\msxfam
\newfam\msyfam
\textfont\msxfam=\tenmsx  \scriptfont\msxfam=\sevenmsx
  \scriptscriptfont\msxfam=\fivemsx
\textfont\msyfam=\tenmsy  \scriptfont\msyfam=\sevenmsy
  \scriptscriptfont\msyfam=\fivemsy
\edef\msx{\hexnumber\msxfam}

\mathchardef\upharpoonright="0\msx16
\let\restriction=\upharpoonright
\def\Bbb#1{\tenmsy\fam\msyfam#1}

\def\restrict{{\restriction}}
\def\re{{\restriction}}

\def\Smallskip{\vskip1.4truecm}
\def\Bigskip{\vskip2.2truecm}

\def\qed{{\vcenter{\hrule height.4pt \hbox{\vrule width.4pt height5pt
 \kern5pt \vrule width.4pt} \hrule height.4pt}}}
\def\notin{{\in}\kern-5.5pt / \kern1pt}
\def\ok{\vbox{\hrule height 8pt width 8pt depth -7.4pt
    \hbox{\vrule width 0.6pt height 7.4pt \kern 7.4pt \vrule width 0.6pt height 7.4pt}
    \hrule height 0.6pt width 8pt}}
\def\nt{{\leq}\kern-1.5pt \vrule height 6.5pt width.8pt depth-0.5pt \kern 1pt}
\def\sd{{\times}\kern-2pt \vrule height 5pt width.6pt depth0pt \kern1pt}
\def\zp#1{{\hochss Y}\kern-3pt$_{#1}$\kern-1pt}

\def\bb{{\frak b}}
\def\cc{{\frak c}}
\def\dd{{\frak d}}
\def\gg{{\frak g}}

\def\hom{{\frak{hom}}}

\def\PPP{{\frak P}}
\def\rr{{\frak r}}

\def\BB{{\Bbb B}}
\def\CC{{\Bbb C}}

\def\GG{{\Bbb G}}
\def\II{{\Bbb I}}
\def\JJ{{\Bbb J}}

\def\LL{{\Bbb L}}

\def\MM{{\Bbb M}}

\def\PP{{\Bbb P}}
\def\QQ{{\Bbb Q}}
\def\RR{{\Bbb R}}
\def\SS{{\Bbb S}}
\def\TT{{\Bbb T}}
\def\UU{{\Bbb U}}
\def\VV{{\Bbb V}}
\def\WW{{\Bbb W}}

\def\sm{{\smallskip}}
\def\ce#1{{\centerline{#1}}}
\def\no{{\noindent}}
\def\la{{\langle}}
\def\ra{{\rangle}}
\def\sub{\subseteq}

\def\ha{{\hat{\;}}}
\def\em{{\emptyset}}
\def\sem{\setminus}
\def\omom{{\omega^\omega}}
\def\omlom{{\omega^{<\omega}}}
\def\omup{{\omega^{\uparrow\omega}}}
\def\omlup{{\omega^{\uparrow < \omega}}}
\def\twoom{{2^\omega}}
\def\twolom{{2^{<\omega}}}
\font\small=cmr8 scaled\magstep0
\font\smalli=cmti8 scaled\magstep0
\font\capit=cmcsc10 scaled\magstep0
\font\capitg=cmcsc10 scaled\magstep1
\font\dunh=cmdunh10 scaled\magstep0
\font\dunhg=cmdunh10 scaled\magstep1
\font\dunhgg=cmdunh10 scaled\magstep2

\font\sanse=cmss10 scaled\magstep0

\font\bolds=cmssdc10 scaled\magstep0

\overfullrule=0pt
\openup1.5\jot

\centerline{}
\Bigskip
\centerline{\dunhgg Strolling through Paradise}
\footnote{}{{\openup-6pt {\small {\smalli 1991
Mathematics subject classification.} 03E05 03E35
03E40 \par {\smalli Key words and phrases.}
$\sigma$--ideal, Marczewski's ideal, nowhere
Ramsey sets, Mycielski's ideal; Sacks forcing,
Miller forcing, Laver forcing, Matet forcing,
Silver forcing, Mathias forcing; covering
number\endgraf}}}
\Bigskip
\centerline{\capitg J\"org Brendle\footnote{$^\star$}
{{\small Supported by the MINERVA--foundation and by
DFG--grant Nr. Br 1420/1--1}}}
\Smallskip
{\baselineskip=0pt {\small \noindent Mathematisches
Institut der Universit\"at, Auf der Morgenstelle 10,
72076 T\"ubingen, Germany}}
\Bigskip
\centerline{\capit Abstract}
\bigskip
\noindent With each of the classical tree--like forcings
adjoining a new real, one can associate a $\sigma$--ideal
on the reals in a natural way. For example, the ideal
of Marczewski null sets $s^0$ corresponds to Sacks forcing
$\SS$, while the ideal of nowhere Ramsey sets $r^0$ corresponds
to Mathias forcing $\RR$. We show (in $ZFC$) that none of
these ideals is included in any of the others. We also
discuss Mycielski's ideal $\PPP_2$, and start an investigation
of the covering numbers of these ideals.
\vfill\eject

{\dunhg Introduction}
\Smallskip
In 1935, E. Marczewski [Mar] introduced on the reals the $\sigma$--ideal
$s^0$, consisting of sets $X \sub 2^\omega$ so that for all perfect
trees $T\sub\twolom$ there is a perfect subtree $S\sub T$ with
$[S] \cap X = \em$, where $[S] := \{ f\in\omom ; \; \forall n \;
(f\re n \in S)\}$ denotes the set of branches through $S$
(see [JMS], [Mi 2], [Ve] and others for recent
results on $s^0$). Similarly a set $X \sub \twoom $ is called {\it
$s$--measurable} iff for all perfect trees $T$ there is a
perfect subtree $S \sub T$ with either $[S] \cap X = \em$ or
$[S] \sub X$. Once forcing was born, the algebra of $s$--measurable
sets modulo $s^0$--sets turned out to be of great interest;
it was first studied by G. Sacks [Sa], and henceforth became
known as {\it Sacks} (or {\it perfect set}) {\it forcing} $\SS$.
Since then, many Sacks--like partial orders have been 
investigated (e.g., Mathias forcing $\RR$, Laver forcing $\LL$,
Miller forcing $\MM$ etc. --- see $\S$ 1 for the definitions),
and it is natural to ask how the corresponding $\sigma$--ideals
(i.e. $r^0$, $\ell^0$, $m^0$, respectively) look like.
We note that the ideal $r^0$ of {\it Ramsey null} (or {\it nowhere
Ramsey}) {\it sets} was first considered by Galvin and Prikry [GP],
and has found a lot of attention over the years (see [AFP], [Br], [Co],
[Ma 2], [Pl] and others), while $\ell^0$ and $m^0$ were looked at
only recently in work of Goldstern, Johnson, Repick\'y, Shelah
and Spinas (see [GJS] and [GRSS]).
\par
One of the fundamental questions one may ask about such ideals is whether
an inclusion relation holds between any two of them. We shall show,
in sections 1 and 2 of the present work, that this is not the
case by constructing in ZFC a set $X \in i^0 \sem j^0$ for each
pair $(i^0 , j^0)$ of such ideals. In case of $(s^0 , r^0)$, this
was done previously under some additional set--theoretic assumptions
by Aniszczyk, Frankiewicz, Plewik, Brown and Corazza (see [AFP],
[Br] and [Co]), and our result answers questions of the latter
[Co, Problems 6 and 10]. In case of $(m^0, r^0)$, this answers a
question of O. Spinas (private communication). --- These results
bear some resemblance to the fact that the ideals corresponding
to Cohen forcing $\CC$ and random forcing $\BB$, the meager sets
${\cal M}$ and the null sets ${\cal N}$, are not included one in the
other. There is even $A\sub \twoom$ with $A\in{\cal M}$ and 
$\twoom \sem A\in {\cal N}$. Two ideals with this property
are called {\it orthogonal}. We also investigate the question which
pairs of the ideals considered in our work are orthogonal and
which are not.
\par
Closely related to these ideals is one of the ideals introduced
by J. Mycielski [My], the $\sigma$--ideal $\PPP_2$, consisting of sets $X
\sub \twoom$ so that for all infinite $A \sub \omega$ the restriction
$X \re A := \{ f \re A ; \; f\in X \}$ is a proper subset of the restriction
$2^A$ of the whole space. $\PPP_2$ is easily seen to be included in the
ideal $v^0$ of {\it Silver null sets} (corresponding to {\it Silver
forcing} $\VV$); we extend the work of section 2 by showing that
it is not included in any of the other previously considered ideals
(Theorem 3.1.). --- Given an ideal ${\cal I}$ on the reals, let 
$cov({\cal I})$ be the size of the smallest ${\cal F} \sub {\cal I}$
covering the reals (i.e. satisfying $\forall f \in \twoom \; \exists 
F \in {\cal F} \; (f\in F)$). We shall prove (Theorem 3.3.) that $v^0$
may be large in comparison with $\PPP_2$ by showing the consistency of
$\omega_1 = cov (v^0) < cov (\PPP_2) = \omega_2 = \cc$. This answers
a question addressed by Cicho\'n, Ros{\l}anowski, Stepr\=ans and
W{\c e}glorz [CRSW, Question 1.3.].
\par
We will conclude our considerations with some remarks concerning
the ideal $r^0_{\cal U}$ of {\it Ramsey null sets} with respect
to a Ramsey ultrafilter ${\cal U}$ in section 4. In particular,
we shall relate the size of the smallest set not in $r^0_{\cal
U}$ to the size of the smallest base of ${\cal U}$ and to
a partition cardinal introduced by Blass [Bl 3, section 6].
\bigskip
{\bolds Notation.} Our set--theoretic notation is fairly standard
(see [Je 1] or [Ku]). $\cc$ denotes the cardinality of the continuum.
Given two sets $A$, $B$, we say that $A$ is almost included in $B$
($A \sub^* B$) iff $A \sem B$ is finite. $\star$ is used for
two--step iterations; we refer to [Bau], [Je 2] and [Sh] for iterated
forcing constructions with countable support.
\par
$\omup$ is the space of strictly increasing functions from $\omega$ to
$\omega$, while $\omlup$ is the set of strictly increasing finite
sequences of natural numbers. For a finite sequence $\sigma$
(i.e. $\sigma \in \twolom, \omlom, \omlup$, or ...), we let $|\sigma|
= dom(\sigma)$, the size (or domain) of $\sigma$, and $rng(\sigma)$,
the range of $\sigma$. $\ha$ is used for concatenation of sequences;
and $\la\ra$ stands for the empty sequence. 
The set of binary sequences of length $n$ is
lexicographically ordered as $\la s_i ; \; i<2^n \ra$  by 
$i <j \iff  s_i (|s|) < s_j (|s|)$ where $s \sub s_i, s_j$.
\par
A tree $T \sub \omlom$ is {\it perfect} iff given $s \in T$, there
are $n \neq m$ and $t \supseteq s$ with $t\ha \la n\ra , t \ha\la m\ra
\in T$; $T$ is {\it superperfect} iff given $s\in T$, there are
$t \supseteq s$ and infinitely many $n$ with $t \ha \la n\ra \in T$.
$split(T) := \{ s \in T ;
\; | \{ n\in\omega ; \; s\ha\la n \ra \in T \} | \geq 2 \}$
is the set of split--nodes of $T$. Given $s \in T$, $succ_T (s)
:= \{ n \in\omega ; \; s\ha\la n \ra \in T \}$ denotes the
immediate successors of $s$ in $T$, while $Succ_T (s) : = \{
t \in T ; \; s \sub t \in split(T) \}$ denotes the successor 
split--nodes of $s$ in $T$. Finally, for $s\in T$, $T_s : = \{
t \in T ; \; s \sub t \;\lor\; t \sub s \}$ is the restriction of
$T$ to $s$.
\bigskip
{\bolds Acknowledgment.} I thank Otmar Spinas for asking me the
question which eventually lead to this work.
\Bigskip 
\vfill\eject

{\dunhg $\S$ 1. Preliminaries}
\Smallskip
{\bolds 1.1.} We will consider the following forcing notions.
\par
\item{---} {\it Sacks forcing} $\SS$ [Je 2, part one, section 3],
also called perfect set forcing:
\par\sm
\ce{$T\in\SS \iff T$ is a perfect tree on $\twolom$}
\par
\ce{$T\leq S \iff T \sub S$}
\sm
\item{---} {\it Miller forcing} $\MM$ [Mi 1] (or superperfect tree
forcing or rational perfect set forcing):
\sm\par
\ce{$T \in\MM \iff T$ is a superperfect tree on $\omlom$}
\par
\ce{$T \leq S \iff T\sub S$}
\sm
\no We note that the conditions in $\MM$ all of whose nodes 
have either infinitely many successor nodes or exactly one
successor node are dense in $\MM$, and henceforth restrict our
attention to such conditions.
\sm
\item{---} {\it Laver forcing} $\LL$ [Je 2, part one, section 3]:
\sm
\par\ce{$T\in\LL \iff T \sub \omlom$ is a tree and $\forall \tau\in T
\; (stem (T) \sub \tau \longrightarrow \exists^\infty n \;
(\tau\ha\la n\ra \in T))$}
\par
\ce{$T\leq S \iff T \sub S$}
\sm
\item{---} {\it Willow tree forcing} $\WW$ (see the end of 1.2. for the
reason for introducing this forcing):
\sm\par
\ce{$ (f,A) \in \WW \iff A \sub [\omega]^{<\omega}$ is infinite
and consists of pairwise disjoint sets $\land$}
\par\ce{$\land\; dom (f) = \omega \sem \cup A \;\land\; ran (f) \sub 2$}
\par
\ce{$ (f,A) \leq (g,B) \iff f \supseteq g \;\land\; \forall a \in A
\; \exists B' \sub B \; (a = \cup B') \;\land$}
\par\ce{$\land\; \forall b\in B
\; (b \sub dom (f) \longrightarrow f\re b $ is constant)
} \sm
\no In this p.o. conditions of the form $(f,A)$ where $A = \{ a_n ;
\; n\in \omega\}$ and $\max (a_n) < \min (a_{n+1})$ are dense,
and we shall always work with such conditions.
\sm
\item{---} {\it Matet forcing} $\TT$ [Ma 1, section 6] (see
also [Bl 1, section 4] and [Bl 2, section 5]):
\sm
\par\ce{$ (s,A) \in \TT \iff s \in \omlup \;\land\; A\sub [\omega]^{
<\omega}$ is infinite $\land\; \forall a\in A \; (\max rng (s)
< \min (a))$}
\par
\ce{ $(s,A) \leq (t,B) \iff s \supseteq t \;\land\; \forall a \in A
\; \exists B' \sub B \; (a = \cup B') \;\land\; \exists B' \sub B
\; (rng (s) \sem rng (t) = \cup B')$}
\sm
\no Again, we may restrict our attention to conditions $(s,A)$
with second coordinate $A = \{ a_n ; \; n\in\omega\}$ satisfying
$\max (a_n) < \min (a_{n+1})$.
\sm
\item{---} {\it Silver forcing} $\VV$ [Je 2, part one, section 3]:
\sm\par
\ce{$f \in \VV \iff dom(f) \sub\omega$ is coinfinite $\land\;
ran(f) \sub 2$}
\par
\ce{$f \leq g \iff f \supseteq g$}
\sm
\item{---} {\it Mathias forcing} $\RR$ [Je 2, part one, section 3]:
\sm\par
\ce{$(s,A) \in \RR \iff s\in\omlup \;\land\; A \sub \omega$
is infinite $\land\; \max rng (s) < \min (A)$}
\par
\ce{$(s,A) \leq (t,B) \iff s \supseteq t \;\land\; A \sub B \;\land
\; rng (s) \sem rng(t) \sub B$}
\sm
\bigskip
{\bolds 1.2.} We note that these forcings are defined on different
underlying sets, e.g. elements of Sacks forcing are subsets of $\twolom$,
elements of Laver forcing subsets of $\omlom$, and Mathias forcing
consists of elements of $\omlup \times [\omega]^\omega$. For our purposes
we need, however, that all these forcings act on the same space,
and we choose $\omup$ to be this space; i.e. we shall think of
each of the forcings as adding a new strictly increasing function
from $\omega$ to $\omega$.
\par
To be more explicit, note first that 
Miller forcing is forcing equivalent to
\sm
\ce{$\MM' = \{ T \sub \omlup ; \; T$ is superperfect$\}$.}
\sm
\no Henceforth, when talking about Miller forcing, we shall mean the
latter p.o. A similar remark applies to Laver forcing.
\par
Next we remark that Mathias and Matet forcing are just uniform
versions of Laver and Miller forcing, respectively, whereas both
$\WW$ and $\VV$ are uniform versions of $\SS$
($\VV$ even being a uniform version of $\WW$). Namely call a Laver
tree $T$ {\it uniform} iff there is $A_T\in [\omega]^{\omega}$
so that for all $\sigma \in T$ extending $stem (T)$, we have 
$succ_T (\sigma) = A_T \sem ( \sigma (|\sigma| -1) +1)$.
Then we can identify $\RR$ and $\RR ' := \{ T\in\LL ;\; T $ is
uniform$\}$: uniform trees $T \in \RR'$ correspond to the pairs
$(stem (T) , A_T) \in \RR$. A similar argument works for the
other forcings.
\par
Finally, define $F : \twolom \to \omlup$ by
\sm
\ce{$ F(\sigma) :=$ the increasing enumeration of $\sigma^{-1} (
\{ 1 \})$.}
\sm
\no $F$ extends to a map $\hat F : 2^\omega \sem \{ f; \; | f^{-1}
(\{ 1 \}) | < \omega \} \to \omup$ defined by
\sm
\ce{$ \hat F (f) := \cup \{ F(f \re n) ; \; n\in\omega\}$.
}\sm
\no $\hat F$ is easily seen to be a homeomorphism. Now, given a
Sacks tree $S \in \SS$, we let 
\sm
\ce{$ \tilde F (S) := \{ F(\sigma) ;\; \sigma \in S \}$.}
\sm
\no Then $\tilde F (S)$ is a perfect subtree of $\omlup$; it is
compact iff $\forall f \in [S] \; ( | f^{-1} (\{ 1\}) |
= \omega)$ [which we can assume, the set of such conditions being dense
in $\SS$]. By a further pruning argument, we may assume all $\tilde F
(S)$ are two--branching; i.e. $\forall \sigma\in\tilde F (S) \;
(| succ_{\tilde F (S)} (\sigma) | \leq 2 )$. Thus the
copy of Sacks forcing on $\omup$ looks exactly like the original
Sacks forcing. 
Henceforth, when talking about Sacks forcing, we shall
mean the p.o. $\{ \tilde F(S) ; \; S\in\SS\}$.
Furthermore we see that every Miller tree is a Sacks
tree (more explicitly, is of the form $\tilde F (S)$ for some
$S\in\SS$).  
Using this (and similar remarks applied to $\VV$ and $\WW$) we see
that the following inclusion relations between the p.o.s
under consideration hold.
\bigskip
\halign{\hskip 6truecm # \quad & \hfill # \quad & \hfill # \quad &
\hfill # \quad & \hfill # \hskip 6truecm \cr
$\SS$ & $\supset$ & $\MM$ & $\supset$ & $\LL$ \cr
\noalign{\sm}
$\cup$ & & $\cup$ &&\cr
\noalign{\sm}
$\WW$ & $\supset$ & $\TT$ & & $\cup$ \cr
\noalign{\sm}
$\cup$ &&& & \cr
\noalign{\sm}
$\VV$ & & $\supset$ && $\RR$ \cr}
\bigskip
\no We realize at this point that $\WW$ arises in a natural way.
The relation between $\SS$ and $\WW$ is like the one between $\MM$
and $\TT$, while the pair $(\SS , \VV)$ corresponds to the pair
$(\LL , \RR)$.
\par
We close this subsection with yet another remark concerning the uniform
forcings. Let $\PP := {\cal P} (\omega) / [\omega]^{<\omega}$,
ordered by almost inclusion, i.e.
\sm
\ce{$A \leq B \iff A \sub^* B$;}
\sm
\no and put $\QQ := {\cal P} ([\omega]^{<\omega}) 
/ [[\omega]^{<\omega}]^{<\omega}$, ordered by
\sm
\ce{$A \leq B \iff \forall^\infty a \in A \; \exists B' \sub B \;
(a = \cup B').$}
\sm
\no Both $\PP$ and $\QQ$ are $\sigma$--closed forcing notions; $\PP$
adjoins a Ramsey ultrafilter ${\cal U}$ on $\omega$ [Ma], while
$\QQ$ adds a stable ordered--union ultrafilter ${\cal V}$ on
$[\omega]^{<\omega}$ (see [Bl 1] for this notion).
It is well--known that $\RR$ is forcing--equivalent to the two--step
iteration $\PP \star \RR_{\breve{\cal U}}$, where $\breve{\cal U}$
is the $\PP$--name for the generic Ramsey ultrafilter, and $\RR_{\cal
U}$ is the $\sigma$--centered Mathias forcing with an ultrafilter
${\cal U}$ [Ma]; similarly $\VV$ decomposes as $\PP \star \GG_{\breve{\cal
U}}$, where $\GG_{\cal U}$ is Grigorieff forcing [Gr].
To this corresponds that $\TT$ is forcing--equivalent to $\QQ \star \TT_{
\breve{\cal V}}$, where ${\cal V}$ is the $\QQ$--name for the generic 
ultrafilter on $[\omega]^{<\omega}$, and $\TT_{\cal V}$
is the $\sigma$--centered Matet forcing with an ordered--union
ultrafilter ${\cal V}$ [Bl 2]; similarly $\WW$ decomposes as $\QQ \star \GG_{\breve
{\cal V}}$, where $\GG_{\cal V}$ is a Grigorieff--like
forcing (we leave the details of this to the reader).
\par
In case of the uniform forcings (i.e. $\RR , \TT , \VV , \WW$),
we sometimes will have to go back to the original notation of the
conditions; we shall always mark the places where we do so, and work
in general with trees.
\bigskip
{\bolds 1.3.} There is a natural way to associate the $\sigma$--ideal
of {\it $\JJ$--null sets} with any of the tree forcings $\JJ$
defined above:
\sm
\ce{$ j^0 := \{ A \sub \omup ; \; \forall T \in \JJ \; \exists
S \leq T \; ([S] \cap A = \em ) \}$.}
\sm
\no Of course, for the forcings with compact trees (i.e. $\SS$,
$\WW$ and $\VV$), one would rather define the corresponding
ideals on $2^\omega$; e.g. Marczewski's ideal $s^0$
(see [Mar] and others) is usually defined as:
\sm
\ce{$ s^0 : = \{ A \sub \twoom ; \; \forall T \in \SS \; \exists
S \leq T \; ( [S] \cap A = \em ) \}$.}
\sm
\no However, putting
\sm
\ce{$\bar F (s^0) := \{ \hat F [A] ; \; A \in s^0 \;\land\;
\forall f \in A \; ( | f^{-1} ( \{ 1 \} ) | = \omega \}$,}
\sm
\no we get the corresponding ideal on $\omup$, and shall 
henceforth work with the latter (and even call it $s^0$).
\par
It is sometimes helpful to think of an $i^0$--set as the complement
of the set of branches of all trees in some dense subset of $\II$
(or some maximal antichain of $\II$).
\par
One of the main goals of this work will be to show that none of these
ideals $i^0$ is included in any other, for the various forcing notions
$\II$ introduced in subsection 1.1. Note that we trivially have $i^0 \sem
j^0 \neq\em$ whenever $\JJ \not\sub \II$: the set of branches of a tree
$T \in \JJ \sem \II$ must be a member of $i^0$, because given any
$S \in\II$, there is $\sigma \in S \sem T$, hence $[S_\sigma ]
\cap [T] = \em$. Thus we are left with showing that $i^0$ is not
included in $j^0$ in case $\JJ \sub \II$. This will be done in
section 2.
\bigskip
{\bolds 1.4.} We make some general remarks concerning the constructions
in section 2. First note that given $\la I_\alpha ; \;
\alpha < \cc \ra \sub \II$ dense and letting $\la J_\alpha ;
\; \alpha < \cc \ra$ be an enumeration of $\JJ$, it suffices to construct
$\la x_\alpha ; \; \alpha < \cc \ra$ so that
\sm
\item{(i)} $x_\alpha \not\in \bigcup_{\beta < \alpha} [I_\beta] \cup
\{ x_\beta ; \; \beta < \alpha \}$;
\par
\item{(ii)} $x_\alpha \in [J_\alpha]$.
\par\sm
\no Then we will have $X = \{ x_\alpha ; \; \alpha < \cc \} \in i^0
\sem j^0$. $X \not\in j^0$ is obvious. To see $X \in i^0$, fix $I \in
\II$; find $\alpha $ so that $I_\alpha \leq I$. Find $\{ I_\beta' ; \;
\beta < \cc \} \sub \II$, an antichain of conditions below $I_\alpha$
with $[I_\beta' ] \cap [I_\gamma'] = \em$ for $\beta \neq \gamma$
(we leave it to the reader to verify that such $I_\beta'$ can be found
for each of our forcings $\II$, using an a.d. family of size $\cc$).
As $| [I_\alpha] \cap X | < \cc$ we necessarily find some $\beta$ with
$[I_\beta'] \cap X = \em$; thus we are done.
\par
In each of our constructions we shall actually see that $| [J_\alpha ]
\sem \bigcup_{\beta < \alpha} [I_\beta ] | = \cc$, so that
the second part of (i) ($x_\alpha \not\in \{ x_\beta ; \;\beta < \alpha \}$)
can be easily satisfied.
\par
The main points of our proofs boil thus down to two steps:
\par
\item{---} choose {\it carefully} a dense set ${\cal I} \sub \II$;
\par
\item{---} find for each $J \in \JJ$ a subtree $J '$ (which will usually be
a homeomorphic copy of $\twolom$) so that $[J' \cap I]$ is small
for all $I \in {\cal I}$ (in general the intersection will be at most
countable).
\par
\no In most cases, it is not difficult to do this. The hardest arguments are
those concerning Marczewski's ideal $s^0$ (in subsections 2.2. and 2.8.).
\par
Finally note that if $\II \supseteq \JJ_1 \supseteq \JJ_2$,
then a set $X \not\in i^0 \sem j^0_1$ constructed along the lines above
will automatically not belong to $j^0_2$ either. Hence we are left
with nine constructions;
they are summarized in the following chart.

\bigskip
{\offinterlineskip \tabskip=0pt
\halign{\strut \vrule# & \quad # \quad & \vrule# & \quad # \quad &
\vrule# & \quad # \quad & \vrule# & \quad # \quad & \vrule# &
\quad # \quad & \vrule# & \quad # \quad & \vrule# & \quad # \quad &
\vrule# & \quad # \quad & \vrule# \cr
\noalign{\hrule}
& $X\in$ && $\ell^0$ && $m^0$ && $s^0$ && $t^0$ && $w^0$ && $r^0$ &&
$v^0$ &\cr
& $\not\in$ &&&&&&&&&&&&&&&\cr
\noalign{\hrule}
& $\ell^0$ &&&& {\bf 2.1.} && (2.2.) && easy && easy && easy && easy &\cr
\noalign{\hrule}
& $m^0$ && easy &&&& {\bf 2.2.} && easy && easy && easy && easy &\cr
\noalign{\hrule}
& $s^0$ && easy && easy &&&& easy && easy && easy && easy &\cr
\noalign{\hrule}
& $t^0$ && easy && {\bf 2.6.} && (2.2./8.) &&&& {\bf 2.3.} && easy && easy
&\cr
\noalign{\hrule}
& $w^0$ && easy && easy && {\bf 2.8.} && easy &&&& easy && easy &\cr
\noalign{\hrule}
& $r^0$ && {\bf 2.5.} && (2.1./6.) && (2.2./8.) && {\bf 2.9.} &&
(2.3./7.) &&&& {\bf 2.4.} &\cr
\noalign{\hrule}
& $v^0$ && easy && easy && (2.8.) && easy && {\bf 2.7.} && easy &&&\cr
\noalign{\hrule} } }
\bigskip
{\bolds 1.5.} Given one of our tree forcings $\JJ$, call a set $A \sub 
\omup$ {\it $j$--measurable} iff for all $T\in\JJ$ there is $S\leq T$
with either $[S] \sub A$ or $[S] \cap A =\em$. If the first alternative
always holds, then $A$ is a {\it $j^1$--set} (this is equivalent to
saying that $\omup \sem A \in j^0$). $A$ is {\it $j$--positive}
iff there is $T\in\JJ$ with $[T] \sub A$. Two ideals ${\cal I}$,
${\cal J}$ on the reals are said to be {\it orthogonal} iff there
is $A\sub\omup$ with $A\in{\cal I}$ and $\omup\sem A\in{\cal J}$.
\par
We note that $\II \supseteq \JJ$ implies that $j$--positive sets
are $i$--positive; in particular $i^0$ and $j^0$ cannot be 
orthogonal. Furthermore, if all sets of reals are $j$--measurable,
then $i^0 \sub j^0$. Using [GRSS, section 2], it is easy to see
that $AD$ implies $\ell$-measurability of all sets of reals; hence it
implies $m^0 \sub \ell^0$ (this observation is due to O. Spinas;
his original argument was somewhat different).
\par
In subsection 2.10. we shall return to the question whether 
$i^0$ and $j^0$ can be orthogonal in case $\II \not\subseteq \JJ$
and $\JJ \not\subseteq\II$.

\Bigskip
\vfill\eject
\input strol2.tex

\Bigskip
\input strol3.tex

\vfill\eject\end

%% file: strol2.tex
%

{\dunhg $\S$ 2. The main results}
\Smallskip
{\bolds 2.1. Theorem.} $m^0 \setminus \ell^0 \neq \emptyset$.
\smallskip
{\capit Proof.} Call a Miller tree $M \in \MM$
an {\it apple tree} iff:
\sm
\ce{$(*_1) \;\;\; \forall \sigma \in split (M) \;
({\rm if} \;\; n>m \;\land\; \sigma\hat{\;}\la n\ra , \sigma\hat{\;}
\la m\ra \in M \;\land\; \sigma\hat{\;}\la m\ra \sub \tau \in
Succ_M (\sigma),$}
\par\ce{$ {\rm then:} \;\;\forall k \in \vert \tau \vert 
\; (\tau (k) < n))$, and}
\sm
\ce{$(*_2) \;\;\; \sigma \subset \tau, \sigma, \tau \in
split (M) \Longrightarrow \vert \tau \vert \geq \vert \sigma \vert +
2$.}
\sm
\no A standard pruning argument shows that given $N \in \MM$
there is an apple tree $M \leq N$. 
\sm
We construct a {\it pear subtree} $P_L = \{ \sigma_t ; \; t \in
2^{<\omega} \} \sub L$ of a Laver tree $L$ which is a copy of
$2^{<\omega}$ such that:
\sm
\item{(I)} $\sigma_{\la\ra} = stem (L)$;
\par
\item{(II)} given $\sigma_t \in L$, $\sigma_{t\hat{\;}
\la 0\ra} = \sigma\ha\la k \ra$ and $\sigma_{t\ha\la 1 \ra}
=\sigma \ha\la \ell\ra$ such that they are in $L$, $\ell > k
> \max \{ \max\; rng (\sigma_{t '} ) ; \; \vert t '\vert
=\vert t\vert \}$.
\par
\sm \no It is immediate from the definition of a Laver tree that
this can be done.
\sm
{\capit Claim.} {\it $\vert [ M \cap P_L ] \vert \leq 1$ whenever
$M$ is an apple tree and $P_L$ is a pear tree.}
\sm
{\it Proof.} Assume to the contrary that $f_1 \neq f_2 \in [ M
\cap P_L]$; fix $\sigma \in \omega^{<\omega}$ so that
$\sigma \sub f_1 , f_2$ and $f_1 (\vert\sigma\vert) < f_2 (\vert
\sigma\vert)$. As both $f_1$ and $f_2$ are branches of $M$,
we must have $f_1 (\vert \sigma\vert + 1) < f_2 (\vert \sigma\vert 
)$; on the other hand, both being branches of $P_L$, we
get $f_1 (\vert\sigma\vert + 1) > f_2 (\vert\sigma\vert)$,
a contradiction. $\qed$
\sm
Now let $\la M_\alpha ;\; \alpha < \cc\ra$ enumerate all
apple trees, and let $\la L_\alpha ; \; \alpha < \cc\ra$ enumerate
all Laver trees. Using the above we easily construct $\la
x_\alpha ;\; \alpha < \cc\ra$ so that 
\par
(i) $x_\alpha \not\in \bigcup_{\beta <\alpha} [M_\beta ] \cup
\{ x_\beta ; \; \beta < \alpha\}$;
\par
(ii) $x_\alpha \in [L_\alpha]$.
\par
\no Then $X = \{ x_\alpha ; \; \alpha < \cc \} \in m^0 \setminus
\ell^0$ by the remarks made in subsection 1.4. $\qed$
\bigskip
{\bolds 2.2. Theorem.} $s^0 \sem m^0 \neq\em$.
\sm
{\capit Proof.} We proceed as before --- but the argument is somewhat
more involved. I.e. we find $\la S_\alpha ;\; \alpha < \cc\ra \sub
\SS$ dense, and construct $\la x_\alpha ; \; \alpha < \cc\ra$
so that
\sm
$x_\alpha \in [M_\alpha]$ and
\par
$x_\alpha \not\in \bigcup_{\beta <\alpha} [S_\beta] \cup
\{ x_\beta ;\;\beta <\alpha\}$,
\sm
\no where $\la M_\alpha ; \;\alpha < \cc\ra$ is an enumeration
of all Miller trees.
\bigskip
\no{\dunh A partition result for Sacks trees}
\sm
We start with thinning out the Sacks trees. Given $S \in \SS$,
$\sigma = stem (S)$ and $i\neq j$ so that $\sigma\ha\la i\ra , \sigma
\ha\la j \ra \in S$, we put $A_S^{\la i,j\ra} : = \{ \la \rho , \tau \ra
; \; \sigma\ha\la i\ra \sub \rho \in split (S) \;\land\; \sigma\ha\la j\ra
\sub \tau \in S \;\land\; |\rho| = |\tau | \}$.
\sm
{\bolds Lemma 1.} {\it Let $S$ be a Sacks tree, $\sigma = stem (S)$,
$i\neq j$ so that $\sigma\ha\la i\ra , \sigma\ha\la j\ra\in S$.
Assume we have a two--place relation $R \sub 
A_S := A_S^{\la i,j\ra}$. Then there is $S' \leq S$ with the
same stem so that \par
\item{---} either $\forall \la\rho,\tau\ra\in A_{S'} \; (\la\rho,\tau\ra\in R)$
\par
\item{---} or $\forall\la\rho,\tau\ra\in A_{S'} \; (\la\rho,\tau\ra\not\in R)$.}
\sm
{\capit Proof.} Assume first:
\sm
\item{$(*)$} there are $n\in\omega$ and $\rho,\tau_k \in S$,
$k < n$, with $\sigma\ha\la i\ra\sub\rho$, $\sigma\ha\la j\ra \sub
\tau_k$ and $\vert\tau_\ell\vert = \vert\tau_k\vert$ such
that for all $\rho ' \supseteq \rho$ in $split (S)$ there is $k<n$
with:
\par\sm
\itemitem{$(\heartsuit)$} $\forall\tau' \supseteq \tau_k \;
((\vert \rho'\vert=\vert\tau'\vert \;\land\; \tau'\in S)
\Longrightarrow \la\rho',\tau'\ra \not\in R)$.
\par\sm
\no Then fix such $n , \rho, \tau_k, k<n$. For $\rho ' \supseteq
\rho$ in $split (S)$ let 
\sm
$k(\rho ') := \min \{ k<n; \; \exists \rho'' \supseteq \rho' $ in
$split(S)$ so that $(\heartsuit)$ holds for $\rho''$ and $k\}$.
\sm
\no Clearly there is $\rho ' \supseteq \rho$ in $split(S)$ so that
$k(\rho '') = k(\rho ') =: k$ for all $\rho '' \supseteq \rho '$
in $split(S)$. Now construct $\la \rho_s ;\; s\in 2^{<\omega}\ra$
so that 
\sm
\item{(i)} $\rho' \sub \rho_s \in split(S)$;
\par
\item{(ii)} $(s \sub t \Rightarrow \rho_s \sub \rho_t)
\;\land\; \rho_{s\ha\la 0\ra} (\vert\rho_s\vert) \neq \rho_{s\ha\la
1\ra} (\vert \rho_s\vert)$;
\par
\item{(iii)} $(\heartsuit)$ holds for $\rho_s$ and $k$.
\par\sm
\no This can be done easily. Let $S' = \{ \rho_s \restrict n ; \;
s\in 2^{<\omega} \;\land\; n\in\omega\} \cup \{\tau\in S ;
\; \tau \sub \tau_k \;\lor\;\tau_k \sub \tau \}$. Clearly the second
alternative of the Lemma holds for $S'$.
\par
So suppose $(*)$ fails; we construct, by recursion on $\vert s\vert$,
$\la \rho_s ;\; s\in 2^{<\omega}\ra$ and  $\la\tau_s ;\; s\in 2^{<\omega}
\ra$ so that
\sm
\item{(a)} $\sigma\ha\la i\ra \sub \rho_s \in split(S) 
\;\land\; \sigma\ha\la j\ra \sub \tau_s \in S$; \par
\item{(b)} $(s \sub t \Rightarrow \rho_s \sub \rho_t , \tau_s \sub
\tau_t ) \;\land\; \rho_{s\ha\la 0\ra} (\vert \rho_s\vert ) \neq
\rho_{s\ha\la 1 \ra} (\vert\rho_s\vert )$ and if $s$ and $t$ are 
incompatible, then so are $\tau_s$ and $\tau_t$;
\par
\item{(c)} $\vert s \vert = \vert t \vert \Rightarrow \vert  \tau_s
\vert=\vert\tau_t\vert \geq\vert\rho_s\vert$ and $\la\rho_s ,\tau_t
\restrict \vert \rho_s\vert\ra\in R$.
\par\sm
\no Assume we are at step $m$ in the construction;
i.e. we have $\la\rho_s ;\; s\in 2^{<m} \ra$, $\la \tau_s ;\; s\in
2^{<m}\ra$ as above. First choose $\la\tilde\tau_t ; \;
t\in 2^m\ra$ and $\la\tilde\rho_t ; \; t\in 2^m\ra \sub S$
pairwise incomparable so that $s\subset t$ implies $\tau_s
\sub \tilde\tau_t$ and $\rho_s \subset \tilde\rho_t$ --- and
also $\tilde\rho_{s\ha\la 0\ra} (\vert \rho_s\vert) \neq
\tilde\rho_{s\ha\la 1\ra} (\vert \rho_s\vert)$ for $s\in 
2^{m-1}$.
Let $\{ t_k ; \; k < 2^m\}$ enumerate $2^m$.
By recursion on $k$ find $\rho_{t_k}$ and $\tilde\tau^k_t$ such that
for all $t\in 2^m$ 
\sm
\item{(A)} $\tilde\rho_{t_k} \sub \rho_{t_k} \in split(S)$;
\par
\item{(B)} $\tilde\tau_t \sub \tilde\tau^{k-1}_t\sub\tilde\tau^k_t$,
$\vert\tilde\tau^k_t\vert =\vert \rho_{t_k} \vert$, and
$\la\rho_{t_k} , \tilde\tau^k_t \ra\in R$.
\par\sm
\no This can be done, because $(*)$ fails for $2^m , \tilde\rho_{t_k} ,
\tilde\tau^{k-1}_t (t\in 2^m)$. Finally put $\tau_t
=\tilde\tau_t^{2^m - 1}$. This completes the
construction.
\par Putting $S' = \{ \rho_s \re n , \tau_s \re n ; \;
s \in \twolom\;\land\; n\in\omega\}$, we see that the first
alternative of the Lemma holds for $S'$. $\qed$
\sm
Now let us assume we have $S\in\SS$ and finitely many pairwise disjoint
relations $R_i \sub\{ \la\rho ,\tau \ra ; \; \rho , \tau \in S \;
\land\; \vert\rho\vert = \vert\tau\vert\}, i<k$, with $\bigcup_{i<k}
R_i = \{\la \rho,\tau\ra ;\;\rho, \tau\in S \;\land\; \vert
\rho \vert=\vert\tau\vert\}$.
We say a splitting node $\sigma\in S$ is {\it of type $\la i,j\ra$}
$(i,j \in k)$ in $S$ iff: letting $n_0  < n_1$ so that $\sigma\ha\la n_0
\ra , \sigma\ha\la n_1\ra\in S$, we have
\sm
$\forall \la\rho, \tau\ra \in A_{S_\sigma}^{\la n_0 , n_ 1 \ra }
\; (\la \rho , \tau \ra \in R_i )$ and
\par
$\forall \la \rho , \tau \ra\in A_{S_\sigma}^{\la n_1 , n_0 \ra} \;
( \la \rho , \tau \ra \in R_j).$
\sm
\no 
Using a standard fusion argument and Lemma 1 we see:
\sm
{\bolds Lemma 2.} {\it Given $S \in \SS$, and $R_i , i<k,$ as above, 
there are $S' \leq S$ and $\la i,j\ra \in k^2$ so that each
splitting node $\sigma \in S'$ is of type $\la i,j\ra$ (in which
case we say $S'$ is of type $\la i,j\ra$).} $\qed$
\sm
\no Given $S\in\SS$ so that $\vert split(S) \cap \omega^n \vert \leq 1$
for all $n\in\omega$ define relations $R_i$, $i<3$, as follows: given
$\sigma \in split (S)$ and $\tau \in S$ with $\vert\sigma\vert = \vert\tau\vert
$ arbitrarily, let $n_0 < n_1$ so that $\sigma\ha\la n_0\ra , \sigma\ha
\la n_1\ra\in S$ and $\tau \sub \tau' \in S$ with $\vert\tau ' \vert
= \vert \tau \vert + 1$ ($\tau '$ being unique), and put
\sm
\ce{$\la \sigma, \tau\ra \in R_0 \Longleftrightarrow n_1 <
\tau ' (\vert\tau\vert)$}
\par
\ce{$\la \sigma, \tau \ra\in R_1 \Longleftrightarrow n_0 \leq
\tau ' (\vert\tau\vert ) \leq n_1$}
\par
\ce{$\la \sigma , \tau \ra\in R_2 \Longleftrightarrow n_0 > \tau '(\vert
\tau\vert)$}
\sm
\no Applying Lemma 2, we get:
\sm
{\bolds Corollary.} {\it The set $\{ S\in\SS ;\; \exists \la i,j \ra \in
3^2 \; (S$ is of type $\la i,j\ra )\}$ is dense in $\SS$. } $\qed$
\bigskip
\no{\dunh Subtrees of Miller trees}
\sm
Assume we are given a family $\Sigma = \la \sigma_s ; \;
s\in 2^{<\omega}\ra\sub\omega^{\uparrow <\omega}$
satisfying
\sm
\item{(I)} $s \subset t \Longrightarrow \sigma_s \subset \sigma_t$;
\sm\par
\item{(II)} $\sigma_{s\ha\la 0\ra} (| \sigma_s | )
< \sigma_{ s \ha\la 1\ra} ( | \sigma_s |)$;
\par\sm 
\item{(III)} given $s\in \twolom , f_i \in\twoom , s\ha\la i\ra\sub f_i
(i\in 2)$, and putting $\phi_i := \bigcup_n \sigma_{f_i \re n}$ we have
for all $n\in\omega$ \par
\ce{$\phi_0 (| \sigma_{f_0 \re (|s| + 2n)}| ) < \phi_1 (|\sigma_{
f_0 \re ( | s| + 2n)} | )$,}
\par
\ce{$\phi_0 (|\sigma_{f_0 \re ( |s| + 2n +1)} | ) > \phi_1 ( | \sigma_{
f_0 \re ( | s| + 2n+1)} | )$,}
\par
\ce{$\phi_1 ( |\sigma_{f_1 \re ( | s|+2n)} | ) > \phi_0 ( | \sigma_{
f_1 \re (|s| + 2n)} |)$,}
\par
\ce{$\phi_1 ( |\sigma_{f_1 \re (|s| +2n+1)}| ) < \phi_0 ( | \sigma_{
f_1 \re ( |s| +2n+1)} | )$;}
\par\sm
\no then we call the closure $C(\Sigma)$ under initial segments
a {\it cherry tree}.
\sm
{\bolds Lemma 3.} {\it A Miller tree $M$ contains a cherry subtree
$C(\Sigma_M)$.}
\sm
{\capit Proof.} We construct by recursion on the levels 
the family $\Sigma_M 
= \la \sigma_s ;\; s \in\twolom \ra\sub split (M)$, so that:
\sm
\item{$(\alpha)$} (I) --- (III) above are satisfied;
\par
\item{$(\beta)$} given $s, t , t' \in \twolom$ with $t(0) = 0,
t' (0) = 1$ and $ | t| = | t'|$ we have:
\par
\ce{$ | \sigma_{s\ha t} | < | \sigma_{s \ha t'} |$ in case $|t|$ is
odd,}
\par
\ce{$|\sigma_{s\ha t} | > | \sigma_{s \ha t'} |$ in case $| t|$
is even.}
\sm
To start, let $\sigma_{\la\ra} : = stem (M)$, and choose splitting
nodes $\sigma_{\la 0\ra} , \sigma_{\la 1\ra} \supseteq \sigma_{\la\ra}$
with $\sigma_{\la 0 \ra} (|\sigma_{\la\ra} | ) < \sigma_{\la 1\ra}
(| \sigma_{\la\ra} | )$ and $| \sigma_{\la 0 \ra} | <
| \sigma_{\la 1 \ra} |$.
\par
Assume $\la \sigma_t ; \; t\in 2^{\leq n} \ra$ have been constructed
satisfying $(\alpha)$ and $(\beta)$ above. Enumerate $\la t_k ;\;
k\in 2^n \ra = 2^n$ in such a way that $k < \ell$ is equivalent to
$| \sigma_{t_k} | > | \sigma_{t_\ell} |$ (this is possible by
$(\beta)$); now recursively find $\sigma_{t_k} \sub \sigma_{t_k \ha
\la i\ra} \in split(M) (i\in 2)$ so that:
\sm
\item{$(\bullet)$} $| \sigma_{t_k \ha\la i \ra } | < | \sigma_{
t_\ell \ha\la j \ra} |$ for $k < \ell$ or $(k=\ell$ and $i<j)$;
\par
\item{$(\bullet)$} $\sigma_{t_k \ha\la 0 \ra} ( | \sigma_{t_k} |
) < \sigma_{t_k \ha\la 1\ra} ( | \sigma_{t_k} | )$;
\par
\item{$(\bullet)$} $\sigma_{t_\ell \ha\la j\ra} ( | \sigma_{t_\ell}
| ) > \sigma_{t_k} ( | \sigma_{t_\ell} | )$ for $k<\ell$
and $j \in 2$.
\par\sm
\no This can be done easily. It is straightforward to verify that
$(\alpha)$ and $(\beta)$ are still satisfied. $\qed$
\sm
Using a similar --- but much easier --- construction, we see:
\sm
{\bolds Lemma 4.} {\it A Miller tree has a subtree of type 
$\la 2 , 2\ra$.} $\qed$
\sm
Unfortunately neither a cherry tree nor a type $\la 2 ,2\ra$--tree
will suffice for our purposes. We have to somehow "amalgamate" these
two types of trees to prove the final lemmata (see below). So
suppose we are given a system $\Sigma = \la \sigma_{\la s,t\ra} ; \; s,
t \in \twolom \;\land\; |s| = |t| \ra \sub \omlup$ such that,
letting $M = M(\Sigma) := \{ \sigma_{\la s,t \ra} \re n ; \;
n\in\omega \;\land\; \sigma_{\la s,t \ra} \in \Sigma \}$ and calling
it a {\it mango tree}, we have:
\sm
\item{\sanse (I)} $(s,t) \subset (s' , t') \Longrightarrow \sigma_{\la s,t\ra}
\subset \sigma_{\la s' , t' \ra}$;
\par
\item{\sanse (II)} $M^f : = \{ \sigma_{\la f\re i , t \ra} \re n ; \; i,n \in\omega
\;\land\; t\in 2^i \}$ is a cherry tree;
\par
\item{\sanse (III)} whenever $f_i , g_i \in \twoom (i\in 2), f_0 \neq f_1 ,
s\sub f_i , f_i (|s|) = i$, then, putting $\phi_i =
\bigcup_n \sigma_{\la f_i \re n, g_i \re n\ra}$, we have
\par\sm
\ce{$\phi_i ( | \sigma_{\la f_i \re n , g_i \re n \ra} | ) >
\phi_j ( | \sigma_{\la f_i \re n , g_i\re n\ra} |)$}
\par\sm
for $i\neq j$ and $n > |s|$, and 
\par\sm
\ce{$\phi_1 ( |\sigma_{\la s , g_1 \re |s|\ra}|) > \phi_0 (|\sigma_{\la
s , g_1 \re |s| \ra } |)$.}
\sm
\no So a mango tree is a kind of "two--dimensional" tree, the vertical sections
of which are cherry trees while the horizontal sections are of
type $\la 2 , 2\ra$ (this is a particular instance of {\sanse (III)},
for $g_0 = g_1$).
\par
To construct a $\Sigma = \la \sigma_{\la s,t\ra} ; \; s,t \in
\twolom \;\land\; |s|=|t|\ra$ giving rise to a mango tree,
proceed as in the proof of Lemma 3, guaranteeing along the way that:
\sm
\item{$(\tilde\alpha)$} {\sanse (I) --- (III)} are satisfied;
\par
\item{$(\tilde\beta)$} given $s,t,t' \in \twolom$ with $t(0) = 0 ,
t' (0) = 1$ and $| t| = | t'|$, and $f\in\twoom$, we have
\par
\ce{$| \sigma_{\la f \re |s| + |t| , s\ha t\ra}| < |
\sigma_{\la f\re |s|+|t| , s\ha t' \ra} |$ in case $|t|$ is odd,}
\par
\ce{$ | \sigma_{\la f\re |s| +|t| , s \ha t\ra} | > |\sigma_{
\la f \re |s| + |t| , s\ha t' \ra} |$ in case $|t|$ is even;}
\par
\item{$(\tilde\gamma)$} in case $s,s' , t, t' \in 2^n$ for some $n$,
and $s$ precedes $s'$ in the lexicographic ordering of $2^n$, we have
\par
\ce{$|\sigma_{\la s,t\ra} | < | \sigma_{\la s' , t' \ra} |$.}
\par\sm
\no In step 0 of the construction, put $\sigma_{\la \ra} : = stem (M)$, and
choose split--nodes $\sigma_{\la i,j \ra} (i,j \in 2)$ extending
$\sigma_{\la\ra}$ with
$\sigma_{\la 0,0 \ra} ( | \sigma_{\la\ra} | ) < \sigma_{\la 0,1\ra}
(\vert \sigma_{\la\ra}| ) < \sigma_{\la 1,0 \ra} (| \sigma_{\la\ra} |)
< \sigma_{\la 1,1 \ra} ( | \sigma_{\la\ra} |)$ and $|\sigma_{\la 0,0\ra }
| < |\sigma_{\la 0,1\ra} | < |\sigma_{\la 1,0\ra} | < |\sigma_{\la 1,1
\ra} |$. --- In step $n$, let $\la s_k ; \; k \in 2^n \ra$ enumerate 
$2^n$ lexicographically; and proceed by recursion on $k$.
For fixed $k$, run the argument in the proof of Lemma 3 twice to get 
$\sigma_{\la s_k \ha\la i \ra , t \ra}$, where $i\in 2, t \in 2^{n+1}$. --- Hence
we proved:
\sm
{\bolds Lemma 5.} {\it A Miller tree contains a mango subtree.} $\qed$
\bigskip
\no{\dunh The final lemmata}
\sm
We are now in a position to conclude our argument by looking at the
intersections of a Sacks tree of one of the types $\la i,j \ra$
$(i,j \in 3)$ with a mango tree. --- Let $E$ denote the set of even
numbers. Given a system $\Sigma = \la \sigma_s ; \; s \in \twolom \ra
\sub \omlup$ satisfying
\sm
\ce{$s \subset t \Longrightarrow \sigma_s \subset \sigma_t$}
\sm
\no (and thus defining a tree $T(\Sigma): = \{ \sigma_s \re n ; \;
s\in \twolom \;\land\; n\in\omega\}$) and a function $f\in 2^E$, we
can form the tree $T (\Sigma_f) := \{ \sigma_s \re n ; \; s \in \twolom \;
\land \; n\in\omega \;\land\; \forall i \in E \cap dom (s) \;
(s(i) = f(i)) \}\sub T(\Sigma)$.
\sm
{\bolds Lemma 6.} {\it Assume $M = M(\Sigma)$ is a mango tree constructed
from the system $\Sigma = \la \sigma_{\la s,t \ra} ; \; s,t \in \twolom
\;\land\; |s|=|t| \ra$, and $S$ is a
Sacks tree of one of the eight types $\la i,j\ra \in 3^2 \sem
\{\la 2,2\ra\}$, then $| \{ f\in\twoom ; \; | [ M^f \cap S ] | \geq
1 \}| \leq \omega$.}
\sm
{\capit Proof.} We look at $\hat T = \{ \la s,t \ra ; \; s,t \in \twolom
\;\land\; |s|=|t| \;\land\; \sigma_{\la s , t\ra} \in M
\cap S \}$.
This is a compact tree in the plane, hence its projection onto the first
coordinate is compact, too, and thus has either at most countably many 
branches or contains a perfect subtree $T$.
In the first case, we are done, so assume the latter.
\par
Put $s := stem (T)$, and note that there must be $t_0 , t_1 \in 2^{|s| + 1}$
so that both $T_i : = \hat T_{\la s \ha\la i \ra , t_i\ra}$ ($i\in 2$)
contain perfect trees. Find incompatible extensions $\la\la s_i^j , t_i^j \ra ;
\; i,j \in 2 \ra$, $\la s_i^j , t_i^j \ra \in T_i$,
and let $\la f_i^j , g_i^j \ra$ be branches of $T_i$ through
$\la s_i^j , t_i^j \ra$. Put (as in {\sanse (III)}) 
$\phi_i^j : = \bigcup_n \sigma_{\la f_i^j \re n , g_i^j \re n \ra}$,
and let $k_i$ be minimal with $\phi^0_i (k_i) \neq \phi^1_i (k_i)$,
It is a consequence of {\sanse (III)} that we must have 
$\phi^j_i (k_i) > \phi_{1-i}^k (k_i)$ for $i,j,k \in 2$. This entails
(by definition of the types) that $S$ is of type $\la 2, 2 \ra$, a
contradiction. $\qed$
\sm
{\bolds Lemma 7.} {\it Assume $S$ is a Sacks tree of type $\la 2 ,2\ra$,
$C = C(\Sigma)$ is a cherry tree constructed from the system
$\Sigma = \la \sigma_s ; \; s\in\twolom \ra$ and $f\in 2^E$, then
$|[ C (\Sigma_f ) \cap S ] | \leq\omega$.}
\sm
{\capit Proof.} Put $C_f = C(\Sigma_f)$ and assume the conclusion
is false. Then $C_f \cap S$ must contain a perfect subtree;
in particular there are $s,t \in\twolom$ so that $\sigma_s , \sigma_t 
\in split (C_f \cap S)$ and $s\ha\la 0 \ra \sub t$. Note that $|s|$ and
$|t|$ must be odd. As $S$ is of type $\la 2 ,2 \ra$ we must have $\sigma_{t
\ha\la i\ra} (| \sigma_t | ) > g ( | \sigma_t| )$ for any 
$g \in [ C_f \cap S ]$ extending $\sigma_{s\ha\la 1 \ra}$. On the other
hand, $C$ being a cherry tree, we have $\sigma_{t\ha\la i\ra} (| \sigma_t
|) < g (|\sigma_t |)$ for any such $g$, a contradiction. $\qed$
\sm
{\bolds Corollary.} {\it If $M$ is a mango tree, and ${\cal S}$ is a family
of less than $\cc$ Sacks trees all of which are of type $\la i,j \ra$
for some $\la i,j \ra \in 3^2$, then $| [M] \sem \bigcup_{S\in{\cal S}}
[S] | = \cc$.}
\sm
{\capit Proof.} First apply Lemma 6 to find $f \in \twoom$ so that
$[M^f \cap S] = \em$ for all trees in ${\cal S}$ which are not
of type $\la 2,2\ra$. Choose $g \in 2^E$ arbitrarily and apply Lemma 7
to find $\cc$ many $\phi \in [C_g ] \sem \bigcup_{s\in{\cal S}}
[S]$, where $C$ is the cherry tree $M^f$. $\qed$
\sm
We can now complete the proof of Theorem 2.2.: let $\la S_\alpha ;
\;\alpha < \cc \ra$ enumerate the Sacks trees of type $\la i,j \ra$ for
some $\la i,j\ra \in 3^2$ --- and construct $\la x_\alpha ; \; \alpha
< \cc \ra$ as required using Lemma 5 and the above Corollary. $\qed$
\bigskip

{\bolds 2.3. Theorem.} $w^0 \sem t^0 \neq \em$.
\sm
{\capit Proof.} This follows from Theorem 3.1. $\qed$
\bigskip

{\bolds 2.4. Theorem.} $v^0 \sem r^0 \neq \em$.
\sm
{\capit Proof.} This follows from Theorem 3.1., too. $\qed$
\bigskip

{\bolds 2.5. Theorem.} $\ell^0 \sem r^0 \neq\em$.
\sm
{\capit Proof.} Call a Laver tree $T \in \LL$ a {\it peach tree}
iff
\sm
\ce{$(\star) \;\;\;$ for any $\sigma, \tau\in T$, if $stem (T)
\sub \sigma, \tau$ and $\sigma \neq \tau$, then $succ_T (\sigma)
\cap succ_T (\tau) = \em$.}
\sm
\no Given $S\in\LL$, there is $T\in\LL$ so that $T \leq S$ and
$T$ is a peach tree (this is a standard fusion argument).
\sm
Construct an {\it orange subtree} $O_M = \{ \sigma_t ; \; t\in 2^{<\omega}
\} \sub M$ of a Mathias tree $M$ as follows:
\sm
\item{(I)} $\sigma_{\la\ra} = stem (M)$;
\par
\item{(II)} suppose $\sigma_t$ for $\vert t \vert \leq n$ is defined;
choose $\ell > k > \max \{ \max rng (\sigma_t) ; \; \vert t \vert
\leq n\}$ such that $k, \ell \in A$ (where $(\sigma_{\la\ra} ,
A)$ is the Mathias condition in usual notation corresponding
to
$M$); then $\sigma_{t\ha\la 0\ra} = \sigma_t \ha\la k\ra$ and 
$\sigma_{t\ha\la 1\ra} = \sigma_t \ha \la\ell\ra$ for any
$t$ with $\vert t \vert = n$.
\par\sm
{\capit Claim.} {\it $\vert [L \cap O_M] \vert \leq 2$ whenever
$L$ is a peach tree and $O_M$ is an orange tree.}
\sm
{\it Proof.} Suppose $f,g,h \in [L \cap O_M]$ were three distinct
elements. Find $n$ such that $f\restrict n \neq g\restrict n \neq
h\restrict n \neq f \restrict n$. Then (without loss) $f (n-1) =
g (n-1)$ and $f \restrict (n-1) \neq g \restrict (n-1)$ --- by the
properties of the orange tree $O_M$; this contradicts the fact
that $L$ is a peach tree. $\qed$
\sm
Using peach and orange trees we complete the proof as in Theorem
2.1. $\qed$
\bigskip

{\bolds 2.6. Theorem.} $m^0 \sem t^0 \neq\em$.
\sm
{\capit Proof.} This is very similar to the proof of Theorem 2.5. A
Miller tree $M$ is a {\it plum tree} iff given $\sigma , \tau \in
split (M)$ distinct, the sets $\bigcup \{ rng (\rho) \sem rng (\sigma)
; \; \rho\in Succ_M (\sigma) \}$ and $\bigcup \{ rng (\rho) \sem
rng (\tau) ; \; \rho \in Succ_M (\tau) \}$ are disjoint. The
set of plum trees is easily seen to be dense in $\MM$.
\sm
Next, given a Matet tree $T$, construct a {\it lemon subtree}
$L_T = \{ \sigma_s \re n ; \; s\in 2^{<\omega} \} \sub T$ so
that
\sm
\item{(I)} $\sigma_{\la\ra} = stem (T)$;
\par
\item{(II)} if $\sigma_s$ for all $s$ of length $\leq n$ are constructed,
choose two finite sets $a , b \sub\omega$ with $\max \{ \max rng (\sigma_s
) ; \; s\in 2^n \} < \min (a) \leq \max (a) < \min (b)$ so that
$a,b \in A_T$, where $A_T$ is the second coordinate in the Matet condition
in usual notation, and put $\sigma_{s\ha\la 0\ra} = \sigma_s \ha \tau_a$,
$\sigma_{s \ha\la 1\ra} = \sigma_s \ha \tau_b$ for all $s\in 2^n$,
where $\tau_a$ ($\tau_b$, resp.) is the increasing enumeration of
$a$ ($b$, resp.).
\par
\sm
\no We see (as in the proof of Theorem 2.5.) that $\vert [L_T \cap M]
\vert \leq 2$ if $M$ is a plum tree and $L_T$ a lemon tree. We conclude
the proof of the Theorem as usual. $\qed$
\bigskip

{\bolds 2.7. Theorem.} $w^0 \sem v^0 \neq \em$.
\sm
{\capit Proof.} Call a willow tree $W$ a {\it fig tree} iff for all
$\sigma \in split (W)$ there are exactly two successor split--nodes
$\tau_1 , \tau_2$ with $\vert \tau_1\vert \geq \vert \tau_2 \vert +2$
and $\tau_1 (\vert \sigma \vert ) < \tau_2 (\vert \sigma \vert )$
(in the language of our original willow conditions this means that
$W$ corresponds to $(f_W , A_W)$ so that, if $A_W =
\{ a_n ; \; n\in\omega\}$ with $\max (a_n) < \min (a_{n+1})$, we
have $\vert a_n \vert \geq 2$ and $\forall n \;\exists i_n\; ( \min
(a_n) < i_n < \min (a_{n+1}) \;\land\; f_W (i_n) =1)$). Clearly
these conditions are dense in $\WW$.
\sm
Now, given a Silver tree $V$, construct a {\it date subpalm}
$D_V$ so that each $\sigma \in split (D_V)$ has exactly two
successor split--nodes $\tau_1 ,\tau_2 $ with $\vert \tau_1 \vert
= \vert \tau_2 \vert$. To do this construct
recursively the split--nodes $\{ \sigma_s ; \; s\in 2^{<\omega} \}$
of $D_V$ as follows:
\sm
\item{(I)} $stem (V) = \sigma_{\la\ra}$;
\par
\item{(II)} assuming $\sigma_s , s\in 2^n$, are constructed, let $i_n , j_n$ be
the $2n$--th and the $(2n +1)$--th elements of $\omega \sem dom (f_V)$
(the corresponding Silver condition in original notation), and put
$\sigma_{s\ha\la 0\ra} = \sigma_s \ha \la i_n \ra\ha \tau \ha\tau'$,
$ \sigma_{s\ha\la 1\ra} = \sigma_s \ha \tau \ha \la j_n \ra\ha\tau'$,
where $\tau, \tau'$ are the increasing enumerations of $f_V^{-1}
(\{ 1\} ) \cap ( i_n , j_n )$, and $f_V^{-1} ( \{ 1 \} )
\cap (j_n , i_{n+1} )$, respectively.
\par\sm
{\capit Claim.} {\it If $W$ is a fig tree and $D_V$ is a date palm,
then $\vert [ W\cap D_V] \vert \leq 1$.}
\sm
{\it Proof.} Assume $f_1 \neq f_2 \in [ W\cap D_V]$. Choose $\sigma$ such
that $\sigma \sub f_1 , f_2$ and $f_1 (\vert \sigma \vert )
< f_2 (\vert\sigma\vert)$. Next let $\tau$ (possibly empty)
be such that $\sigma \ha\la f_1 (\vert\sigma\vert ) \ra \ha \tau \sub
f_1$ and $\sigma \ha\tau \sub f_2$ and $f_1 (\vert \sigma \vert +
\vert \tau \vert + 1) \neq f_2 (\vert \sigma \vert + \vert\tau\vert )$.
$f_1 , f_2 \in [W]$ requires that $f_1 (\vert \sigma \vert + \vert \tau
\vert + 1 ) < f_2 ( \vert\sigma \vert + \vert \tau \vert )$, while
$f_1 , f_2 \in [D_V]$ requires that $f_1 ( \vert\sigma \vert + \vert \tau
\vert + 1 ) > f_2 (\vert\sigma\vert + \vert\tau\vert)$. $\qed$
\sm
We complete the proof of the Theorem as usual. $\qed$
\bigskip

{\bolds 2.8. Theorem.} $s^0 \sem w^0 \neq \em$.
\sm
{\capit Proof.} We shall use the same dense subset of $\SS$ as in subsection
2.2., namely the set $\{ S\in\SS ; \; \exists \la i,j \ra \in 3^2 \;
( S$ is of type $\la i,j \ra ) \}$.
\bigskip
\no{\dunh Subtrees of willow trees of type $\la i,j \ra$}
\sm
{\bolds Lemma 1.} {\it A willow tree $W$ has a type $\la 0,0 \ra$--subtree
$T_W$.}
\sm
{\capit Proof.} Assume $W$ corresponds to the pair $(f_W ,A_W)$ in the usual
notation for willow conditions. Without loss $A_W = \{ a_n ; \; n
\in\omega\}$, $\max (a_n) < \min (a_{n+1} )$, $2 \cdot | a_n | < | a_{n+1}|$,
and also $f_W^{-1} ( \{ 1 \} ) \cap ( \min (a_n) , \min (a_{n+1} )) \neq
\em$ (otherwise go over to a stronger condition).
\par
Let $\la I_n ; \; 1\leq n\in\omega\ra$ be a partition of $\omega$ into intervals
of size $2^{n+1}, \max (I_n) + 1 = \min (I_{n+1})$.
$\tau^n_j  \; (j < 2^{n+1})$ is the increasing enumeration of
the set $a_i \cup [ f_W^{-1} (\{ 1\}) \cap (\min (a_i) , \min (a_{i+1} ))]$
and $\rho_j^n  \; (j < 2^{n+1} )$ is the increasing
enumeration of $f_W^{-1} ( \{ 1 \} ) \cap (\min (a_i) , \min
(a_{i+1}))$,
where $i$ is the $j$--th element of $I_n$. Define recursively $\la \sigma_s
; \; s \in \twolom \ra \sub W$:
\sm
\ce{$\sigma_{\la\ra} = stem (W)$;}
\par
\ce{$\sigma_{\la 0 \ra} = \sigma_{\la \ra} \ha \tau_0^1$; and}
\par
\ce{$\sigma_{\la 1 \ra} = \sigma_{\la\ra} \ha \rho_0^1 \ha \rho_1^1 \ha 
\tau_2^1$.}
\sm
\no Assume $\la \sigma_s ; \; s \in 2^{\leq n} \ra$ have been defined,
$n \geq 1$. Let $\la s_i ; \; i < 2^n \ra$ be the lexicographic 
enumeration of $2^n$.
We put 
\sm
\ce{$\sigma_{s_i \ha\la 0\ra} = \sigma_{s_i} \ha \tau_{2i + 1}^n \ha
\rho_{2i+2}^n \ha ... \ha \rho^n_{2^{n+1}-1} \ha \rho^{n+1}_0 \ha
... \ha \rho_{4i - 1}^{n+1} \ha \tau_{4i}^{n+1}$;}
\sm
\no and
\sm
\ce{$\sigma_{s_i \ha \la 1\ra} = \sigma_{s_i} \ha \rho_{2i+1}^n \ha
... \ha \rho_{2^{n+1} - 1}^n \ha \rho_0^{n+1} \ha ... \ha \rho_{4i+1}^{n+1} \ha 
\tau_{4i + 2}^{n+1}$.}
\sm
Set $T_W := \{ \sigma_s \re n ; \; n\in\omega \;\land\; s\in \twolom\}$.
Note that the construction was set up in such a way that whenever $\sigma
\in T_W$ is a splitnode, then $\sigma$ is of the form $\sigma_s$ for
some $s \in \twolom$; thus the final part of this sequence is some $\tau^n_j$.
By our requirements on the $| a_m |$, this entails
that $\sigma_s$ is longer that any other sequence in the tree
ending in the corresponding $\rho_j^n$. Hence, if $\sigma_1 ,
\sigma_2$ are two immediate successors of $\sigma_s$, then
$\sigma_1 ( | \sigma_s | ) , \sigma_2 ( | \sigma_s | ) < \tau (
| \sigma_s | )$ for any sequence $\tau$ in the tree $T_W$ which is incomparable
with $\sigma_s$. Therefore $T_W$ is of type $\la 0,0 \ra$. $\qed$
\sm
We leave to the reader the proof of the following --- easier ---
result:
\sm
{\bolds Lemma 2.} {\it A willow tree $W$ has a type $\la 0,2 \ra$--subtree $S_W
$ --- in fact, we can construct a subtree $S_W$ of $W$ with
splitnodes $\la \sigma_s ; \; s\in \twolom \ra$ so that $s \subset t$ implies
$\sigma_s \subset \sigma_t$ which satisfies: whenever $t_0 \supseteq s
\ha \la 0\ra$ and $t_1 \supseteq s \ha \la 1\ra$, then $\forall n
\geq | \sigma_s | \; ( \sigma_{t_0} (n) < \sigma_{t_1} (n))$.}
$\qed$
\sm
As in subsection 2.2. we want a kind of "two--dimensional" subtree of
a willow tree so that the sections in the two directions are of type
$\la 0,0\ra$ and of type $\la 0,2\ra$, respectively. To this
end we construct a system $\Sigma = \la\sigma_{\la s,t \ra} ; \;
s,t \in\twolom \;\land\; |s| = |t| \ra \sub W$ so that,
letting $P = P(\Sigma) := \{ \sigma_{\la s,t\ra} \re n ; \; n\in
\omega \;\land\; \sigma_{\la s,t\ra} \in \Sigma \}$ and calling
it a {\it poplar tree}, we have:
\sm
\item{(I)} $s' \subset s , t' \subset t \Longrightarrow \sigma_{\la s'
, t '\ra} \subset \sigma_{\la s , t \ra }$;
\sm
\par
\item{(II)} $P^f := \{ \sigma_{\la f\re i , t \ra} \re n ; \; i,n \in
\omega \;\land\; t\in 2^i \}$ is of type $\la 0,0 \ra$;
\sm\par
\item{(III)} $P_g := \{ \sigma_{\la s , g \re i \ra} \re n ; \; i,n\in\omega
\;\land\; s \in 2^i \}$ is as in Lemma 2 (in particular is
of type $\la 0,2\ra$);
\sm\par
\item{(IV)} whenever $f_i , g_i \in \twoom (i\in 2), f_0 \neq f_1,
s\sub f_i , f_0 (|s| ) =0$, and $f_1 (|s| ) = 1$, then,
putting $\phi_i = \bigcup_n \sigma_{\la f_i \re n , g_i \re n\ra}$,
we have
\par
\ce{$\phi_0 (m) < \phi_1 (m)$ for $m \geq | \sigma_{\la f_1 \re |s| + 1 ,
g_1 \re |s| +1 \ra}|$.}
\par\sm
\no ((IV) is a kind of strengthening of (III) which we shall need in
Lemma 4 below; in fact, (III) itself won't be used.) This construction
is done in a similar fashion as the construction in the proof of
Lemma 1. Namely, we make the same initial assumptions about 
$W = (f_W , A_W)$. Next we take $\la I_s ; \; s\in
\twolom \;\land\; |s|\geq 1 \ra$ a partition of $\omega$ into
intervals of size $2^{|s| +1}$ so that $\max (I_s) +1 \leq
\min (I_t)$ whenever $|s| <|t|$ and $\max (I_s) + 1 \leq
\min (I_t)$ for $s,t$ with $|s|=|t|$ iff $s$ precedes $t$ in the
lexicographic ordering of $2^{|s|}$. We let $\la\tau_j^s  ; \; 
j < 2^{|s| + 1} \ra
$ and $\la\rho_j^s ;\; j < 2^{|s| + 1}\ra$ be defined accordingly.
---
Construct recursively $\la \sigma_{\la s , t \ra} ;\; s,t \in 2^{<\omega}
\;\land\; |s| = |t| \ra \sub W$:
\sm
\ce{$\sigma_{\la\ra} = stem (W)$,}
\par
\ce{$\sigma_{\la 0,0 \ra}  = \sigma_{\la\ra} \ha \tau_0^{\la 0 \ra}$,}
\par
\ce{$\sigma_{\la 0,1 \ra} = \sigma_{\la\ra} \ha \rho_0^{\la 0 \ra}
\ha \rho_1^{\la 0 \ra} \ha \tau_2^{\la 0 \ra}$,}
\par
\ce{$\sigma_{\la 1,0 \ra } = \sigma_{\la\ra} \ha \tau_0^{\la 0 \ra}
\ha \rho_1^{\la 0\ra}
\ha ... \ha \rho_3^{\la 0 \ra} \ha \tau_0^{\la 1 \ra}$,}
\par
\ce{$\sigma_{\la 1,1\ra} = \sigma_{\la\ra} \ha \rho_0^{\la 0 \ra}
\ha ... \ha \rho_1^{\la 1 \ra} \ha \tau_2^{\la 1 \ra}$.}
\sm
\no Assume $\la \sigma_{ \la s,t \ra} ; \; s,t \in 2^{\leq n}
\;\land\; |s| = |t| \ra$ have been defined, $n \geq 1$.
Let $\la s_i ; \; i < 2^n \ra$ be the lexicographic enumeration
of $2^n$. We put:
\sm
\ce{$\sigma_{\la s_i , s_j \ra \ha \la 0,0 \ra} = \sigma_{\la s_i
, s_j \ra } \ha \tau_{2j+1}^{s_i} \ha \rho_{2j + 2}^{s_i} \ha ...
\ha \rho_{2^{n+1} -1}^{s_i} \ha \tau_0^{s_{i+1}} \ha ... \ha \tau_{
2^{n+1} - 1}^{s_{2^n - 1}} \ha \rho_0^{s_0 \ha \la 0 \ra} \ha ...
\ha \rho_{4j - 1}^{s_i \ha\la 0 \ra} \ha \tau_{4j}^{s_i \ha\la 0 \ra}$,}
\par
\ce{$\sigma_{\la s_i , s_j \ra \ha \la 0,1 \ra } = \sigma_{\la s_i
, s_j \ra} \ha \rho_{2j + 1}^{s_i} \ha ... \ha \rho_{2^{n+1} -1}^{s_i}
\ha \tau_0^{s_{i+1}} \ha ... \ha \tau_{2^{n+1} - 1}^{s_{2^n - 1}} \ha
\rho_0^{s_0 \ha \la 0 \ra} \ha ... \ha \rho_{4j + 1}^{s_i \ha\la 0 \ra}
\ha \tau_{4j + 2}^{s_i \ha\la 0 \ra}$,}
\par
\ce{$\sigma_{\la s_i , s_j \ra \ha \la 1, 0 \ra} = \sigma_{\la s_i ,
s_j \ra} \ha \tau_{2j+1}^{s_i}\ha \rho_{2j + 2}^{s_i} \ha ...
\ha \rho_{2^{n+1} -1}^{s_i} \ha \tau_0^{s_{i+1}} \ha ... \ha \tau_{
2^{n+1} - 1}^{s_{2^n - 1}} \ha \rho_0^{s_0 \ha \la 0 \ra} \ha ...
\ha \rho_{4j - 1}^{s_i \ha\la 1 \ra} \ha \tau_{4j}^{s_i \ha\la 1 \ra}$,}
\par
\ce{$\sigma_{\la s_i , s_j \ra \ha \la 1, 1 \ra} = \sigma_{\la s_i ,
s_j \ra}\ha \rho_{2j + 1}^{s_i} \ha ...
\ha \rho_{2^{n+1} -1}^{s_i} \ha \tau_0^{s_{i+1}} \ha ... \ha \tau_{
2^{n+1} - 1}^{s_{2^n - 1}} \ha \rho_0^{s_0 \ha \la 0 \ra} \ha ...
\ha \rho_{4j + 1}^{s_i \ha\la 1 \ra}
\ha \tau_{4j + 2}^{s_i \ha\la 1 \ra}$.}
\sm
\no As in the proof of Lemma 1, we verify that this $\Sigma =
\la \sigma_{\la s , t \ra } ; \; s , t \in \twolom \;\land\; |s|=|t| \ra$
satisfies (I) --- (IV):
\sm
{\bolds Lemma 3.} {\it A willow tree $W$ has a poplar subtree 
$P(\Sigma_W)$.} $\qed$
\sm
We can now complete the proof of Theorem 2.8. along similar lines
as the proof of 2.2. with the following lemmata:
\sm
{\bolds Lemma 4.} {\it Assume $P = P(\Sigma_W) $ is a poplar subtree
of a willow tree $W$, and $S$ is a Sacks tree of one of the seven
types $\la i,j \ra \in 3^2 \sem \{ \la 0,2 \ra , \la 2 , 0 \ra \}$.
Then $| \{ f \in \twoom ; \; | [ S \cap P^f ] | \geq 1 \} |
\leq\omega$.}
\sm
{\capit Proof.} We look at $\hat T :=\{\la s,t \ra ; \; s,t \in\twolom \;\land\;
|s|=|t| \;\land\; \sigma_{\la s ,t \ra} \in S \cap P \}$.
This is a compact tree in the plane, hence its projection onto the
first cooradinate is compact, too, and thus has either at most
countably many branches or contains a perfect subtree $T$. In the
first case, we are done, so assume the latter.
\par
Put $s : = stem (T)$, and note that there must be $t_0 , t_1 \in
2^{|s| + 1}$ so that both $T_i := \hat T_{\la s \ha \la i \ra , t_i \ra}
$ ($i\in 2$) contain perfect trees. Find incompatible extensions $\la\la s_i^j ,
t_i^j \ra ; \; i,j < 2 \ra$, $\la s_i^j , t_i^j \ra \in T_i$,
and let $\la f_i^j , g_i^j \ra$ be branches of $T_i$ through $\la
s_i^j , t_i^j \ra$. Put (as in (IV)) $\phi_i^j : = \bigcup_n
\sigma_{\la f_i^j \re n , g_i^j \re n \ra}$, and let $k_i$ be minimal
with $\phi^0_i (k_i) \neq \phi^1_i (k_i)$. It is a consequence of (IV)
that we must have $\phi_0^j (k_i) < \phi_1^k (k_i)$ for $i,j,k\in 2$.
This entails (by definition of the types) that $S$ is either of type
$\la 0,2 \ra$ or of type $\la 2 , 0 \ra$, a contradiction. $\qed$
\sm
The argument of the following result is similar, but much easier
(just note that $P^f$ is of type $\la 0,0 \ra$ (II)):
\sm
{\bolds Lemma 5.} {\it Assume $P=P(\Sigma_W)$ is a poplar subtree of a 
willow tree $W$, and $S$ is a Sacks tree of type $\la 0,2 \ra$
or $\la 2,0 \ra$. Then for all $f\in\twoom$, $| [S\cap P^f ] |
\leq\omega$.} $\qed$
\sm
Putting Lemmata 4 and 5 together we get:
\sm
{\bolds Corollary.} {\it Assume $P=P(\Sigma_W)$ is a poplar
subtree of a willow tree $W$, and ${\cal S} \sub \{ S \in \SS ;
\; S$ is of type $\la i,j \ra$ for some $\la i,j \ra \in 3^2 \}$
is a family of Sacks trees of size $< \cc$, then 
$|[P] \sem \bigcup_{S \in {\cal S}} [S]|=\cc$.}
\sm
{\capit Proof.} Use Lemma 4 to find $f\in\twoom$ so that $[ S\cap P^f ]
= \em$ for all $S\in {\cal S}$ of one of the types $\la i,j \ra \in
3^2 \sem \{ \la 0,2 \ra , \la 2,0 \ra \}$. Then use Lemma 5
to find many $g$'s so that $\phi = \bigcup_n \sigma_{\la f\re n , 
g \re n \ra}$ is as required. $\qed$
\sm
Now we can conclude the proof of Theorem 2.8. with
the usual argument. $\qed$
\bigskip

{\bolds 2.9. Theorem.} $t^0 \sem r^0 \neq \em$.
\sm
{\capit Proof.} A Matet tree $T \in \TT$ is a {\it maple tree}
iff given $\sigma, \tau \in T$, $\sigma (n) = \tau (n)$
implies $\sigma \re n = \tau \re n$.
\sm
{\capit Claim 1.} {\it Given $T\in \TT$, there is a maple tree $S \leq T$.
}
\sm
{\it Proof.} Let $(s, A_T)$ be the Matet condition in usual notation
associated with $T$; i.e. $s\in \omlup$ and $A_T = \{ a_n ; \;
n \in \omega \}$ is an infinite
set of finite subsets of $\omega$ satisfying $\max rng (s) < \min
(a_0) < ... < \max (a_n) < \min (a_{n+1}) < ...$ for $n\in\omega$.
Going over to a stronger condition, if necessary, we may assume that
whenever $\sigma , \tau\in\omlup$ are distinct, then
$$ \sum_{i < |\sigma|} |a_{\sigma (i)} | \neq \sum_{i<|\tau|}
|a_{\tau (i)} |.$$
[to do this it suffices to guarantee that 
$|a_n| > \sum_{m<n} |a_m|$ for all $n$]. Note that this entails
$|split (T) \cap \omega^n | \leq 1$ for all $n\in\omega$.
\par
To see that such a $T$ must be a maple tree, take $\sigma, \tau \in T$
and $n\in\omega$ with $\sigma\re n \neq \tau\re n$, choose $m_\sigma ,
m_\tau \leq n$ maximal so that $\sigma\re m_\sigma , \tau\re m_\tau
\in split (T)$. In case $m_\sigma = m_\tau$, we necessarily have $\sigma
\re m_\sigma = \tau \re m_\sigma$, and $\sigma (n) \in a_i$,
$\tau(n) \in a_j$ for $i\neq j$, and the values must be distinct.
In case $m_\sigma < m_\tau$ (without loss), $\sigma (n)$ is the
$(n-m_\sigma)$--th value of some $a_i$, and $\tau (n)$ is the
$(n - m_\tau)$--th value of some $a_j$, and again the values
must be distinct. $\qed$
\sm
Now recall the definition of an orange tree from the proof of
Theorem 2.5., and note that the following claim is proved
exactly as the corresponding claim there.
\sm
{\capit Claim 2.} {\it If $T$ is a maple tree and $O_M$ is an orange
tree, then $| [T  \cap O_M] | \leq 2$.} $\qed$
\sm
We conclude the proof of the Theorem as usual. $\qed$   
\bigskip

{\bolds 2.10.} We conclude this section with a result and a question
about orthogonality of our ideals.
\sm
{\bolds Proposition.} {\it The pairs of ideals $(m^0 , w^0)$, $(m^0 , v^0)$,
$(\ell^0 , w^0)$, $(\ell^0 , v^0)$ and $(v^0, t^0)$ are orthogonal.}
\sm
{\capit Proof.} The first four follow from Proposition 3.2. To see the
last, call $V\in\VV$ an {\it oak tree} iff given $i_0 < i_1$, both
in $\omega\sem dom(f_V)$, there is $i_2 \in f_V^{-1} (\{ 1\})$ with
$i_0 < i_2 < i_1$ (here $f_V$ denotes the Silver condition in usual
notation associated with $V$). $T\in\TT$ is an {\it almond tree} iff
$|a_n| \geq 2$ and $\max (a_n) < \min (a_{n+1})$, where $(s, \{ a_n ;
\; n\in\omega\})$ is the Matet condition in usual notation associated
with $T$. The oak trees are dense in $\VV$ and the almond trees are
dense in $\TT$. It is easily seen that $|[T\cap V]| \leq 1$ in case $V$ is
an oak tree and $T$ is an almond tree. Hence, if $\la V_\alpha ;\;\alpha < \cc\ra
=\{$oak trees$\}$ and $\la T_\alpha ; \;\alpha < \cc\ra =\{$almond trees$\}$,
we can easily construct $\la V_\alpha ' ; \;\alpha <\cc\ra$ and $\la T_\alpha
' ; \; \alpha <\cc\ra$ so that for all $\alpha, \beta$:
\sm
\item{(i)} $T_\alpha ' \leq T_\alpha$, $V_\alpha ' \leq V_\alpha$; \par
\item{(ii)} $[T_\alpha ' \cap V_\beta ] = \em$. \par
\sm
\no Then $\bigcup_{\alpha <\cc} [T_\alpha]$ is a $t^1$--set in $v^0$.
$\qed$
\sm
{\bolds Question.} {\it Are $\ell^0$ and $t^0$ orthogonal?}
\sm
\no As neither $\LL \sub \TT$ nor $\TT \sub \LL$, a positive answer
seems plausible.

%% file: strol3.tex
%

{\dunhg $\S$ 3. Some results concerning the Mycielski ideal $\PPP_2$}
\Smallskip
{\bolds 3.1.} Given $A \sub \twoom$ and $X \in [\omega]^\omega$
we let $A \re X := \{ f \re X ; \; f \in A \}$. Recall that the
{\it Mycielski ideal $\PPP_2$} is defined as follows:
\sm
\ce{$\PPP_2 := \{ A \sub \twoom ; \; \forall X \in [\omega]^\omega \;
(A \re X \neq 2^X ) \}$.}
\sm
\no $\PPP_2$ is easily seen to be a $\sigma$--ideal on the reals 
which is contained in both $v^0$ and $w^0$ (see [CRSW] for more
on $\PPP_2$).
\par
As usual we shall be concerned with the isomorphic copy of $\PPP_2$
in the space $\omup$,
\sm
\ce{$ \bar F (\PPP_2) := \{ \hat F [A] ; \; A \in \PPP_2 \;\land\;
\forall f \in A \; ( | f^{-1} (\{ 1 \}) | = \omega ) \}$}
\sm
\no (compare subsection 1.3.), and actually mean the latter ideal
when talking about $\PPP_2$. 
\par
It follows from the results in section 2 (and $\PPP_2 \sub v^0 , w^0$)
that $i^0 \sem \PPP_2 \neq
\em$ for any of the ideals $i^0$ considered so far;
similarly $\PPP_2 \sem s^0 ( m^0 , \ell^0 )\neq \em$ are easy to see.
We conclude this cycle of results by showing:
\sm
{\bolds Theorem.} {\it $\PPP_2 \sem t^0 \neq \em$ (and thus $\PPP_2 \sem
r^0 \neq \em$).}
\sm
{\capit Note.} Theorems 2.3. and 2.4. immediately follow from this
result.
\sm
{\capit Proof.} Given $X \in [\omega]^\omega$ and $ f\in 2^X$ with
$|f^{-1} ( \{ 1 \} ) | =\omega$, let $S_f := \{ g\in 2^\omega ; \;
f \sub g \}$ and $R_f := \{ \hat F (g) ; \; g\in S_f \} \sub \omup$,
the {\it closed set associated with $f$} (see 1.2. for the map $\hat F$).
Let $\la X_\alpha ; \; \alpha < \cc \ra$ enumerate $[\omega]^\omega$
and $\la T_\alpha ; \; \alpha < \cc\ra = \TT$. We shall construct
$\la f_\alpha ; \; \alpha < \cc \ra$ and $\la y_\alpha ; \; \alpha
< \cc\ra$ so that for all $\alpha$:
\sm
\item{(i)} $f_\alpha \in 2^{X_\alpha}$, $|f_\alpha^{-1} ( \{ 1 \} ) | =
\omega$;
\par
\item{(ii)} $y_\alpha \in [T_\alpha ]$;
\par
\item{(iii)} $y_\alpha \not\in \bigcup_{\beta < \cc} R_{f_\beta} 
\cup \{ y_\beta ; \; \beta < \alpha \}$.
\par
\sm
\no Then we clearly have $Y: = \{ y_\alpha ; \; \alpha < \cc \}
\not\in t^0 , r^0$. But also $Y \cap R_{f_\alpha} = \em$ for
all $\alpha < \cc$, and thus $Y \in \PPP_2$.
Suppose we are at step $\alpha$ in the construction;
then we can easily find $f_\alpha \in 2^{X_\alpha}$ so that
$|f_\alpha^{-1} ( \{ 1 \} ) | = \omega$
and $\hat F^{-1} (y_\beta) \re X_\alpha \neq f_\alpha$ for all
$\beta < \alpha$ (this entails $y_\beta \not\in R_{f_\alpha}$).
Thus the main point is to find $y_\alpha \in [T_\alpha] \sem
(\bigcup_{\beta \leq \alpha} R_{f_\beta} \cup
\{ y_\beta ; \; \beta < \alpha \})$.
\par
Given a Matet tree $T \in\TT$ construct an {\it elm subtree}
$E_T = \{ \sigma_s \re n ; \; s \in \twolom \;\land\; n \in
\omega\}$ as follows:
\sm
\item{(I)} $\sigma_{\la\ra} = stem (T)$;
\par
\item{(II)} if $\sigma_s $ are constructed for all $s$ of length $\leq 
n$, choose $2^{n+1}$ distinct finite subsets of $\omega$,
$\la a_i ; \; i \in 2^{n+1} \ra$ with $\max \{ \max rng (\sigma_s ) ;
\; s \in 2^n \} < \min (a_i)$ and $a_i \in A_T$ ($i\in 2^{n+1}$) ---
where $A_T$ is the second coordinate in the Matet condition in usual
notation ---, and put
\sm\par
\ce{$\sigma_{s_i \ha\la 0 \ra} = \sigma_{s_i} \ha \tau_{2i} , \;
\sigma_{s_i \ha\la 1 \ra} = \sigma_{s_i} \ha \tau_{2i + 1}$,}
\sm\par
\item{} where $\{ s_i ; \; i \in 2^n \}$ enumerates $2^n$ and $\tau_j$
is the increasing enumeration of $a_j$ ($j\in 2^{n+1}$).
\par\sm
{\capit Claim.} {\it $| [E_T] \cap R_f | \leq 1$, whenever $E_T$ is
an elm tree and $f\in 2^X$  ($X\in[\omega]^\omega$) is such that
$|f^{-1} ( \{ 1 \} ) | = \omega$.}
\sm
{\it Proof.} Simply note that two distinct branches of $E_T$ have
almost disjoint range while $f^{-1} ( \{ 1 \} )$ is contained in the 
range of any real in $R_f$. $\qed$
\sm
Using this claim we easily conclude the construction of $y_\alpha$.
$\qed$
\bigskip
{\bolds 3.2. Proposition.} {\it $\PPP_2$ is orthogonal to $m^0$ and
$\ell^0$; it is not orthogonal to any of the other ideals.}
\sm
{\capit Proof.} Note that if $i^0$ and $\PPP_2$ are orthogonal,
then $i^0$ and $w^0$ ($v^0$) are orthogonal as well. This 
proves non--orthogonalitiy for most ideals.
\par
We proceed to show the orthogonality of $(\PPP_2 , m^0)$
(the orthogonality of $(\PPP_2 , \ell^0)$ is proved in a similar
fashion). We use the notation in the proof of Theorem 3.1.
Let $\la X_\alpha ;\;\alpha <\cc\ra$ enumerate $[\omega]^\omega$
and $\la M_\alpha ; \; \alpha < \cc\ra = \{ M\in\MM , \;
M$ is a plum tree as in 2.6.$\}$. We construct recursively $\la f_\alpha ;
\; \alpha < \cc\ra$ and $\la M_\alpha ' ; \; \alpha < \cc\ra$
so that for all $\alpha$:
\sm
\item{(i)} $f_\alpha \in 2^{X_\alpha}$, $|f_\alpha^{-1} (\{
1\}| =\omega$; \par
\item{(ii)} $M_\alpha ' \leq M_\alpha$; \par
\item{(iii)} $[M_\alpha '] \cap \bigcup_{\beta < \cc} R_{f_\beta}
=\em$. \par
\sm
\no Then $Y:= \bigcup_{\alpha < \cc} [M_\alpha ']$ is an $m^1$--set
which lies in $\PPP_2$. Assume we are at step $\alpha$; find
$f_\alpha \in 2^{X_\alpha}$ so that $|f_\alpha^{-1} (\{ 1 \} )
| =\omega$ and $f_\alpha \notin \hat F^{-1} ([M_\beta ']) \re
X_\alpha$ for $\beta < \alpha$ (this is easy to do because for
$f\neq g \in [M_\beta ']$, $|(\hat F^{-1} (f))^{-1} (\{ 1\})
\cap (\hat F^{-1} (g))^{-1} (\{ 1\}) | <\omega$). Next note
that, as in the proof of 3.1., we have $|[M_\alpha] \cap R_{f_\beta}
| \leq 1$ for $\beta \leq \alpha$. Thus we find $M_\alpha '
\leq M_\alpha$ with $[M_\alpha '] \cap R_{f_\beta} = \em$
for $\beta \leq\alpha$. This concludes the construction. $\qed$
\bigskip

{\bolds 3.3.} We next turn our attention to consistency results
by showing that $v^0$ may be "much larger" than $\PPP_2$:
\sm
{\bolds Theorem.} {\it It is consistent with $ZFC$ that $\omega_1
= cov (v^0) < cov (\PPP_2) = \omega_2 = \cc$.}
\sm
This Theorem will be proved with an iterated forcing construction
with countable support of length $\omega_2$ over a model for $CH$.
We start with the definition of the forcing we want to iterate,
{\it Sacks forcing with uniform levels} $\UU$:
\sm
\ce{$ S \in \UU \iff S \in \SS \; \land\; \exists A_S \in [\omega]^\omega
\; \forall \sigma \in S \; ( \sigma \in split (S) \longleftrightarrow
|\sigma| \in A_S )$}
\par
\ce{$S \leq T \iff S \sub T$}
\sm
\no (note that we work in $\twolom$ and $\twoom$ in this subsection).
Standard arguments show that $\UU$ is proper (satisfies even axiom
A) and $\omega^\omega$--bounding (see [Bau], [Je 2] or [Sh] for these
notions). We shall be interested in showing that it shares with
Sacks forcing $\SS$ yet
another nice property preserved in countable support iterations
(which is not shared by Silver forcing).
To do this recall that a non--principal ultrafilter ${\cal U}$ on $\omega$ is
called {\it $P$--point} iff for all $\la A_n ; \; n \in \omega \ra$
with $A_n \in {\cal U}$ there is $B \in {\cal U}$ with $B \sub^*
A_n$ for all $n$. A forcing notion $\PP \in V$ is called 
{\it $P$--point--preserving} iff whenever ${\cal U} \in V$ 
is a $P$--point, then ${\cal U}$ generates a $P$--point in $V[G]$,
where $G$ is $\PP$--generic over $V$. 
\sm
{\bolds Lemma 1.} {\it $\UU$ is $P$--point--preserving.}
\sm
{\capit Proof.} Let ${\cal U}$ be a $P$--point.
It is well--known that it follows from the properness
of $\UU$ that it suffices to show that ${\cal U}$ generates an
ultrafilter in $V[G]$ (see, e.g., [BlSh, Lemma 3.2.]).
\par
Choose $S \in \UU$, and let $\dot B$ be a $\UU$--name for an infinite
subset of $\omega$. Let $\dot\chi$ be the $\UU$--name for
the characteristic function of $\dot B$. Let $\la n_i ; \; i\in\omega\ra$
be the increasing enumeration of $A_S$. We may assume that if $\sigma \in
S \cap 2^{n_i}$, then $S_\sigma$ decides $\dot\chi \re i$ (otherwise go over
to a stronger $S' \leq S$, using a standard fusion argument).
\par
With each $f\in [S]$ we associate a set $B_f \sub \omega$
which reflects $f$'s opinion about $\dot B$:
\sm
\ce{$ i \in B_f \iff S_{f \re n_{i+1}} \forces \dot\chi
(i) = 1$}
\sm
\no (notice that $S_{f\re n_{i+1}}$ decides $\dot\chi (i)$). Assume
that
\sm
\ce{$(*) \;\;\;\;\; $ for all $\sigma \in S$ there is $f\in [S]$
so that $B_f \in {\cal U}$.}
\sm
\no Then construct a decreasing sequence of sets $\la B_i \in {\cal U} ; \;
i\in\omega\ra$ as follows:
\sm
\item{$\bullet$} choose $f\in [S]$ so that $B_f \in {\cal U}$,
and put $B_0 := B_f$;
\par
\item{$\bullet$} assume $B_{i-1} $ ($i \geq 1$) is constructed;
find for each of the $2^i$ splitting nodes $\sigma $ of $S$ of length
$2^{n_i}$ an $f_\sigma \in [S]$ extending $\sigma$ with $B_{f_\sigma}
\in {\cal U}$, and let $B_i$ be the intersection of $B_{i-1}$ and these
$2^i$ $B_f$'s.
\par\sm
\no In case $(*)$ fails for some $\sigma \in S$, we can make a similar
construction of sets $\la C_i \in {\cal U} ; \; i \in\omega\ra$ below
$\sigma$, this time taking complements of sets of the form $B_f$, and
intersecting them. The rest of the proof is the same,
whether or not $(*)$ holds, so assume the former is the case.
\par
${\cal U}$ being a $P$--point, we find $B \in {\cal U}$ which is
almost included in all $B_i$. Construct a sequence $\la k_i ;
\; i\in\omega\ra$ so that $k_0 = 0$ and $B \sem B_{k_i + 1}
\sub k_{i+1}$. Note that either $B' := B \cap \bigcup_i [ k_{2i} ,
k_{2i + 1} ) \in {\cal U}$ or $B \cap \bigcup_i [k_{2i + 1} ,
k_{2i + 2} ) \in {\cal U}$. Without loss assume the
former holds. We construct $S' \leq S$ with $A_{S'} = \{ n_{k_{2i+1}} ;
\; 
i\in\omega \}$ such that whenever $f \in [S']$, then $B' \sem
[k_0 , k_1 ) \sub B_f$ (and thus $S' \forces B' \sub^* \dot B$):
\sm
\item{$\bullet$} let $\sigma_{\la\ra} \in S$ of length $n_{k_1}$ be
arbitrary; choose two split--nodes $\tilde\sigma_{\la 0 \ra}
, \tilde \sigma_{\la 1\ra}$ of $S$ extending $\sigma_{\la\ra}$
of length $n_{k_1 + 1}$; let $f_{\tilde\sigma_{\la 0 \ra}}$ and
$f_{\tilde\sigma_{\la 1 \ra }}$ be as in the above construction
(thus we will have $B_{f_{\tilde\sigma_{\la j\ra}}} \supseteq B_{k_1
+1} \supseteq B' \sem [k_0 , k_1)$), and put $\sigma_{\la j\ra} 
:= f_{\tilde\sigma_{\la j \ra}} \re n_{k_3}$;
\par
\item{$\bullet$} assume $\la \sigma_s ; \; s \in 2^i \ra$, $i \geq 1$,
are constructed, and of length $n_{k_{2i + 1}}$; choose two split--nodes
$\tilde\sigma_{s\ha\la j \ra}$ ($j\in 2$) of $S$ extending $\sigma_s$ of
length $n_{k_{2i + 1} + 1}$; let $f_j \in [S]$ be the extension of
$\tilde\sigma_{s\ha\la j\ra}$ in the above construction (thus we
will have $B_{f_j} \supseteq B_{k_{2i + 1} + 1} \supseteq B' \sem
[k_0 , k_1 )$), and put $\sigma_{s\ha\la j \ra} := 
f_j \re n_{k_{2i + 3}}$.
\par\sm
\no This completes the construction, and thus the proof of the Lemma.
$\qed$
\sm
{\bolds Lemma 2.} (Shelah; see, e.g., [BlSh, Theorem 4.1.])
{\it If $\la \PP_\alpha , \dot\QQ_\alpha ; \; \alpha < \lambda \ra$
is a countable support iteration of proper forcing notions so
that
\sm
\ce{$\forces_{\PP_\alpha} "\dot\QQ_\alpha$ is $P$--point--preserving$"$,}
\sm
\no then $\PP_\lambda$ is $P$--point--preserving as well.} $\qed$
\sm
The {\it reaping number} $\rr$ is the size of the smallest family ${\cal F}$
of infinite subsets of $\omega$ so that given any $A \in [\omega]^\omega$,
there is $B \in {\cal F}$ with either $B \sub A$ or $ B \cap A 
=\em$. A translation to cardinal invariants
of the well--known fact that Silver forcing
$\VV$ adjoins a new subset of $\omega$ which neither contains nor
is disjoint from an old infinite subset of $\omega$ leads to:
\sm
{\bolds Lemma 3.} (Folklore) {\it $cov (v^0) \leq \rr$.}
\sm
{\capit Proof.} Define $G: \twolom \to \twolom$ as follows:
\sm
\ce{$G (\la 0 \ra) = \la 0 \ra$}
\par
\ce{$G (\la 1 \ra) = \la 1 \ra$}
\par
$$G(\sigma) = \cases{ G(\sigma \re (|\sigma| - 1)) \ha \la \sigma (| \sigma
| -2) \ra & if $\sigma(|\sigma| - 1 ) = 0$ \cr
G(\sigma \re (| \sigma | - 1 )) \ha\la 1-\sigma (|\sigma| - 2 ) \ra &
if $\sigma( |\sigma| -1 ) = 1$ \cr}$$
Let $\hat G : \twoom \to \twoom$ be defined by:
\sm
\ce{$\hat G (f) := \cup \{ G (f \re n) ; \; n\in\omega\}$.}
\sm
\no Then $\hat G$ is a homeomorphism of $2^\omega$. Now let
${\cal F} \sub [\omega]^\omega$ be a witness for $\rr$ of size $\rr$.
Given $A\in[\omega]^\omega$, let $\tilde A_i = 
\{ f\in \twoom ; \; \forall n\in A \; (f(n) = i) \}$,
and let $\hat A_i = \hat G^{-1} [\tilde A_i]$ ($i=0,1$).
\par
We claim that $\hat A_i \in v^0$. To see this let $f\in \VV$; choose
$n\in\omega$ minimal with $n\not\in dom(f)$. Let $m > n$ be minimal
with $m \in A$. There are (at least) two ways to extend $f$ to
$\hat f$ so that $(m+1) \sub dom (\hat f)$; according to which way we choose
we either have $\hat G (g) (m) =0 $ for all $g \supseteq \hat f$
or $\hat G (g) (m) = 1$ for all $g \supseteq \hat f$. Choosing the
first we get $\hat G (\{ g ; \; g \supseteq \hat f \} ) \cap \tilde A_1
= \em$, i.e. $\{ g; \; g\supseteq \hat f \} \cap \hat A_1 = \em$. Thus
$\hat A_1 \in v^0$. Similarly for $\hat A_0$.
\par
Next, $\{ \tilde A_i ; \; A \in {\cal F}\;\land\; i\in 2\} 
$ is
covering; hence so is the inverse image under $\hat G$. $\qed$
\sm
{\capit Proof of the Theorem.} Our model is the generic extension
by the countable support
iteration of $\UU$ of length $\omega_2$ over a model
for $CH$. By properness, no cardinals
are collapsed and $\cc = \omega_2$ in the resulting model.
By Lemmata 1 and 2, there is a $P$--point generated by $\omega_1$
sets in the latter model; the base of this $P$--point is a reaping
family, and thus witnesses $\rr = \omega_1$; hence --- by Lemma 3 ---
$cov(v^0) = \omega_1$.
\par
Next note that, letting $u^0$ be the $\sigma$--ideal corresponding to
$\UU$ (as in 1.3.), we have $\PPP_2 \sub u^0$ (this comes from
the uniform levels); thus $cov (\PPP_2) \geq cov (u^0)$. However,
a standard L\"owenheim--Skolem argument shows that $cov (u^0)
= \omega_2$ in our model
(see [JMS, Theorem 1.2.] for the corresponding
result for Sacks forcing $\SS$ and Marczewski's ideal $s^0$);
hence $cov(\PPP_2) =\omega_2$, too. 
$\qed$
\bigskip

{\bolds 3.4.} We have seen in 3.3. that for two specific 
ideals ${\cal I}$, ${\cal J}$ with ${\cal I} \not\subseteq
{\cal J}$, $cov ({\cal I}) < cov ({\cal J})$ is consistent.
We conjecture that this is true for any such combination of
our ideals. In many cases this is quite easy to see. For example, to
show the consistency of $cov(\ell^0) < cov(m^0)$, simply 
iterate Miller forcing $\omega_2$ times. Then a L\"owenheim--Skolem
argument shows $cov(m^0) =\omega_2$ in the final model;
on the other hand, it is well--known [Mi 1] that $\bb
=\omega_1$ holds, where $\bb$, the {\it unbounding number},
is the smallest family ${\cal F}$ of functions from
$\omega$ to $\omega$ so that for every $g\in\omom$ there
is $f\in{\cal F}$ which is infinitely often above $g$.
An easy translation (like in Lemma 3 in 3.3.) of the folklore
fact that Laver forcing adds a dominating real shows
$cov(\ell^0) \leq \bb$; hence $cov(\ell^0) =\omega_1$.
However, the consistency of $cov(m^0) < cov(\ell^0)$
is open --- and so are several similar problems.
\par
As $w^0 \sub\PPP_2$, the result in 3.3. leads to the following,
more specific, problem.
\sm
{\bolds Question.} {\it Is $cov (w^0) =\omega_1$ in the model
of Theorem 3.3.? Is $cov (w^0) < cov(\PPP_2)$ consistent?}
\sm
\no A positive answer seems plausible.
\Bigskip

{\dunhg $\S$ 4. Odds and ends}
\Smallskip
{\bolds 4.1.} Recall that a non--principal ultrafilter ${\cal U}$ on $\omega$
is {\it Ramsey} (or {\it selective}) iff given a partition
$\pi : [\omega]^2 \to 2$ there is $A \in {\cal U}$ {\it homogeneous
for $\pi$} (this means $| \pi [A]^2 | = 1 $) iff given a partition
of $\omega$ into infinitely many (infinite) pieces $\la A_n ; \;
n\in\omega\ra$ not in ${\cal U}$ there is $A\in{\cal U}$ with
$| A \cap A_n | \leq 1$ for all $n$ (see [Je 1, Lemma 38.1] for the
latter equivalence). Ramsey ultrafilters may not exist [Je 1,
Theorem 91]; but for the remainder of this section we assume
there is at least one Ramsey ultrafilter, called ${\cal U}$.
$\RR_{\cal U}$ is {\it Mathias forcing with respect to ${\cal U}$};
i.e. $(s,A) \in\RR_{\cal U}$ iff $(s,A) \in \RR$ and $A\in{\cal U}$,
the ordering being inherited from $\RR$ (as usual we can think
of $\RR_{\cal U}$ as forcing with certain subtrees of $\omlup$).
With $r_{\cal U}^0$ we denote the corresponding $\sigma$--ideal;
i.e.
\sm
\ce{$X \in r^0_{\cal U} \iff  X \sub \omup \;\land\;
\forall T \in \RR_{\cal U} \;\exists S \leq T \; (S \in \RR_{\cal U}
\;\land\; X \cap [S] = \em )$.}
\sm
\no $r^0_{\cal U}$ was studied by Louveau [Lo] and Corazza [Co].
It is easy to see that $r_{\cal U}^0 \sem i^0 \neq \em$ for any of the
ideals $i^0$ considered so far; on the other hand, it follows from
our results in section 2 that $i^0 \sem r^0_{\cal U} \neq \em$, as well
(in case $i^0 = s^0$, this answers [Co, Problem 7]).
\par
Given a $\sigma$--ideal ${\cal I}$ on the reals, let $non ({\cal I})$
be the size of the smallest set not in ${\cal I}$ (note that
$non ({\cal I})$ is sometimes denoted by $unif ({\cal I})$, the
{\it uniformity} of ${\cal I}$). It is well--known that $non ({\cal I})
= \cc$ for any of the ideals ${\cal I}$ investigated in sections 2 and 3;
Corazza [Co, Theorem 4.3.] proved that $non (r^0_{\cal U}) = \omega_1 < \cc$
is consistent. We shall try to relate $non (r^0_{\cal U})$
to more familiar cardinal invariants. --- Given an ultrafilter
${\cal V}$ on $\omega$, let $\gg ({\cal V})$, the {\it generating
number} of ${\cal V}$, be the size of the smallest base for ${\cal V}$.
The {\it homogeneity number} $\hom$ [Bl 3, section 6] is the size
of the smallest family ${\cal H}$ of subsets of $\omega$ so that
for each partition $\pi : [\omega]^2 \to 2$ there is $H \in {\cal H}$
homogeneous for $\pi$.
\sm
{\bolds Proposition.} {\it $\hom \leq non (r^0_{\cal U}) \leq \gg ({\cal
U})$.}
\sm
{\capit Proof.} Assume $\{ A_\alpha ; \; \alpha < \gg ({\cal U}) \}$
generates ${\cal U}$. Let $\chi_\alpha \in \twoom$ be the characteristic
function of $A_\alpha$. Put $X := \{ f\in \twoom ; \;
\exists \alpha < \gg ({\cal U}) \; \forall^\infty n \;
(f(n) = \chi_\alpha (n) ) \} \sub 2^\omega$.
Clearly $| X | = \gg ({\cal U})$; on the other hand $X \not\in
r_{\cal U}^0$ (here, we think of $r_{\cal U}^0$ as an ideal on
$\twoom$ instead of $\omom$, i.e. we really work with $\bar F^{-1}
(r^0_{\cal U} )$ --- see 1.3.). To see this take $(s,A) \in \RR_{\cal
U}$ arbitrarily; as $A\in{\cal U}$, choose $\alpha < \gg ({\cal U})$
with $A_\alpha \sub A$; let $f\in\twoom$ be such that $f(n) = 1$ iff
$n \in rng(s) \cup A_\alpha$. Then $f\in X$, but $f$ also lies in the subset
of $\twoom$ defined by the condition $(s,A)$ (by this we mean
$\hat F (f) \in T$, where $T$ is the subtree of $\omlup$ corresponding
to the condition $(s,A)$). This proves $non (r_{\cal U}^0 ) \leq
\gg ({\cal U})$.
\par
To see $\hom \leq non (r^0_{\cal U} )$, let $\kappa < \hom$ and
$X = \{ x_\alpha ; \; \alpha < \kappa \} \sub\twoom$.
Put $A_\alpha := x_\alpha^{-1} ( \{ 1 \} )$. Choose a partition
$\pi : [\omega]^2 \to 2$ so that no $A_\alpha \sem n$ is homogeneous
for $\pi$ ($\alpha < \kappa , n\in\omega$). ${\cal U}$ being Ramsey,
there is $A\in{\cal U}$ homogeneous for $\pi$. Now note that given
$(s,B) \in \RR_{\cal U}$ arbitrarily, $[T] \cap \hat F [X] =
\em$ for the tree $T$ associated with the condition $(s, A\cap B)$.
Hence $X \in r_{\cal U}^0$. $\qed$
\sm
{\bolds Corollary.} (MA) {\it $non (r_{\cal U}^0) = \cc$.}
$\qed$
\sm
(This answers [Co, Problem 9].)
\sm
{\bolds Question.} {\it Does either $non (r_{\cal U}^0) = \gg ({\cal U})$
or $non (r_{\cal U}^0) = \hom$ always hold? Or are there models $V_1 ,
V_2$ of $ZFC$ containing Ramsey ultrafilters ${\cal U}_1 \in V_1$
and ${\cal U}_2 \in V_2$ so that $V_1 \models non (r_{{\cal U}_1}^0) <
\gg ({\cal U}_1)$ and $V_2 \models \hom < non (r_{{\cal U}_2}^0)$?}
\bigskip
{\bolds 4.2.} We shall comment on the relation between $\hom$ and
other cardinal invariants of the continuum (though this is not directly
related to our principal object of study, certain $\sigma$--ideals on the
reals, we feel that what follows rounds off our work).
--- Let $\rr_\sigma$ be the size of the smallest family ${\cal F}$
of subsets of $\omega$ so that for each countable sequence $\la
Y_n ; \; n\in\omega\ra$ of infinite subsets of $\omega$ there is $F \in
{\cal F}$ almost contained in either $Y_n$ or $\omega\sem Y_n$
for all $n$. Recall the definition of $\rr$ from subsection 3.3.
Finally let $\dd$ be the {\it dominating number}, the size of the 
smallest family ${\cal F}$ of functions from $\omega$ to $\omega$ so
that every function $f\in\omom$ is dominated everywhere by a 
function in ${\cal F}$. Blass proved [Bl 3, Theorem 17]
\sm
\ce{$ \max (\rr , \dd ) \leq \hom \leq \max ( \rr_\sigma , \dd ).$}
\sm
\no We show (answering [Bl 3, Question 5] positively):
\sm
{\bolds Proposition.} {\it $\hom = \max ( \rr_\sigma , \dd )$.}
\sm
{\capit Proof.} By Blass' result it suffices to prove that $\rr_\sigma
\leq \hom$. To do this, let ${\cal H} \sub [\omega]^\omega$
be such that for all $\pi : [\omega]^2 \to 2$
there is $H \in {\cal H}$ homogeneous for $\pi$ $(\star)$.
It suffices to show:
\sm
\ce{$(\star\star) \;\;\;\;\; \forall \la Y_n ;\; n\in\omega\ra
\sub [\omega]^\omega \; \exists X \in{\cal H} \; \forall n\in\omega \;
(X\sub^* Y_n \;\lor\; X \sub^* \omega \sem Y_n )$.}
\sm
\no To do this we associate with $\la Y_n ;\; n\in\omega\ra$ a partition
$\pi : [\omega]^2 \to 2$ as follows.
\sm
{\sanse step 0.} We first define $\pi$ on all pairs $\{ x,y\}$ so that
$x \in Y_0$ and $y \in \omega \sem Y_0$. For such pairs let
$$\pi (\{ x,y \} ) = \cases{0 & iff $x<y$ \cr
1 & iff $x>y$\cr}$$
\par
{\sanse step 1.} Next we define $\pi$ on all pairs $\{ x,y\}$ so that
\par\item{(i)} $\pi $ has not yet been defined on $\{ x,y\}$;
\par\item{(ii)} $x\in Y_1$ and $y \in \omega \sem Y_1$.
\par
\no For such pairs let
$$\pi (\{ x,y \} ) = \cases{0 & iff $x<y$ \cr 1 & iff $x>y$ \cr}$$
etc...
\par
If, in the end, $\pi$ has not yet been defined on all pairs $\{ x,y \}$
we define it arbitrarily on the remaining pairs. By $(\star)$
choose $H\in{\cal H}$ homogeneous for $\pi$. Let us check that
$H$ satisfies condition $(\star\star)$. To reach a contradiction,
assume there is $n\in\omega$ so that
\sm
\ce{$ | H \cap Y_n | =\omega = |H \cap (\omega\sem Y_n) | \;\;\; (+)$.}
\sm
\no Let $n_0$ be the least such $n$. Let $k\in\omega$ be so that
for all $n < n_0$ either $H \sem k \sub Y_n $ or $H \sem k \sub \omega
\sem Y_n$ $(++)$. By $(+)$ there are $x_0 , x_1 , x_2 \in H$
so that $k \leq x_0 < x_1 < x_2$ and $x_0 \in Y_{n_0}$,
$x_1 \in \omega \sem Y_{n_0}$ and $x_2 \in Y_{n_0}$. By $(++)$
$\pi$ was not defined for the pairs $\{ x_0 , x_1\}$ and $\{ x_1 ,
x_2 \}$ in the steps up to step $n_0 -1$ of the construction of
$\pi$. Hence in step $n_0$ of the construction we set $\pi ( \{
x_0 , x_1 \} ) = 0$ and $\pi ( \{ x_1 , x_2 \}) = 1$.
This contradicts the homogeneity of $H$. $\qed$
\Bigskip

{\dunhg References}
\Smallskip
\itemitem{[AFP]} {\capit B. Aniszczyk, R. Frankiewicz and S. Plewik,}
{\it Remarks on (s) and Ramsey--measurable functions,} Bulletin
of the Polish Academy of Sciences, Mathematics, vol. 35 (1987), pp. 
479-485.
\smallskip
\itemitem{[Bau]} {\capit J. Baumgartner,} {\it Iterated forcing,}
Surveys in set theory (edited by A.R.D. Mathias), Cambridge
University Press, Cambridge, 1983, pp. 1-59.
\smallskip
\itemitem{[Bl 1]} {\capit A. Blass,} {\it Ultrafilters related to Hindman's
finite unions theorem and its extensions,} in: Logic and Combinatorics,
ed. S. Simpson, Contemporary Mathematics, vol. 65 (1987), pp. 89-124.
\smallskip
\itemitem{[Bl 2]} {\capit A. Blass,} {\it Applications of superperfect
forcing and its relatives,} in: Set theory and its applications,
ed. J. Stepr\=ans and S. Watson, Springer Lecture Notes in Mathematics,
vol. 1401 (1989), pp. 18-40.
\smallskip
\itemitem{[Bl 3]} {\capit A. Blass,} {\it Simple cardinal
characteristics of the continuum,} in: Set Theory of the Reals
(H. Judah, ed.),
Proceedings of the Bar--Ilan conference in honour of
Abraham Fraenkel, 1993, pp. 63-90.
\smallskip
\itemitem{[BlSh]} {\capit A. Blass and S. Shelah,}
{\it There may be simple $P_{\aleph_1}$-- and $P_{\aleph_2}$--points
and the Rudin--Keisler ordering may be downward directed,}
Annals of Pure and Applied Logic, vol. 33 (1987),
pp. 213-243.
\smallskip
\itemitem{[Br]} {\capit J. Brown,} {\it The Ramsey sets and related
sigma algebras and sets,} Fundamenta Mathematicae, vol. 136 (1990),
pp. 179-185.
\smallskip
\itemitem{[CRSW]} {\capit J. Cicho\'n, A. Ros{\l}anowski, J.
Stepr\=ans and B. W{\c e}glorz,} {\it Combinatorial properties of
the ideal $\PPP_2$,} Journal of Symbolic Logic, vol. 58 (1993),
pp. 42-54.
\smallskip
\itemitem{[Co]} {\capit P. Corazza,} {\it Ramsey sets, the 
Ramsey ideal, and other classes over $\RR$,} Journal of
Symbolic Logic, vol. 57 (1992), pp. 1441-1468.
\smallskip
\itemitem{[GP]} {\capit F. Galvin and K. Prikry,} {\it Borel sets
and Ramsey's theorem,} Journal of Symbolic Logic, vol. 38 (1973),
pp. 193-198.
\smallskip
\itemitem{[GJS]} {\capit M. Goldstern, M. Johnson and O. Spinas,}
{\it Towers on trees,} to appear in Proceedings of the American
Mathematical Society.
\smallskip
\itemitem{[GRSS]} {\capit M. Goldstern, M. Repick\'y, S.
Shelah and O. Spinas,} {\it On tree ideals,} preprint.
\smallskip
\itemitem{[Gr]} {\capit S. Grigorieff,} {\it Combinatorics on ideals
and forcing,} Annals of Mathematical Logic, vol. 3 (1971), pp. 
363-394.
\smallskip
\itemitem{[Je 1]} {\capit T. Jech,} {\it Set theory,} Academic Press,
San Diego, 1978.
\smallskip
\itemitem{[Je 2]} {\capit T. Jech,} {\it Multiple forcing,}
Cambridge University Press, Cambridge, 1986.
\smallskip
\itemitem{[JMS]} {\capit H. Judah, A. Miller and S. Shelah,}
{\it Sacks forcing, Laver forcing, and Martin's Axiom,}
Archive for Mathematical Logic, vol. 31 (1992), pp. 145-162.
\smallskip
\itemitem{[Ku]} {\capit K. Kunen,} {\it Set theory,} North-Holland,
Amsterdam, 1980.
\smallskip
\itemitem{[Lo]} {\capit A. Louveau,} {\it Une m\'ethode topologique
pour l'\'etude de la propri\'et\'e de Ramsey,} Israel Journal of
Mathematics, vol. 23 (1976), pp. 97-116.
\smallskip
\itemitem{[Mar]} {\capit E. Marczewski,} {\it Sur une classe
de fonctions de W. Sierpi\'nski et la classe correspondante
d'ensembles,} Fundamenta Mathematicae, vol. 24 (1935), pp. 17-34.
\smallskip
\itemitem{[Ma 1]} {\capit P. Matet,} {\it Some filters of partitions,}
Journal of Symbolic Logic, vol. 53 (1988), pp. 540-553.
\smallskip
\itemitem{[Ma 2]} {\capit P. Matet,} {\it Happy families and completely
Ramsey sets,} Archive for Mathematical Logic, vol. 32 (1993), pp.
151-171.
\smallskip
\itemitem{[Ma]} {\capit A. R. D. Mathias,} {\it Happy families,}
Annals of Mathematical Logic, vol. 12 (1977), pp. 59-111.
\smallskip
\itemitem{[Mi 1]} {\capit A. Miller,} {\it Rational perfect
set forcing,} Contemporary Mathematics, vol. 31 (Axiomatic
Set Theory, 1984, edited by J. Baumgartner, D. Martin and
S. Shelah), pp. 143-159.
\smallskip
\itemitem{[Mi 2]} {\capit A. Miller,} {\it Special subsets of the
real line,} Handbook of set--theoretic topology, K. Kunen and 
J. E. Vaughan (editors), North--Holland, Amsterdam, 1984, pp.
201-233.
\smallskip
\itemitem{[My]} {\capit J. Mycielski,} {\it Some new ideals of subsets
on the real line,} Colloquium Mathematicum, vol. 20 (1969), pp. 71-76.
\smallskip
\itemitem{[Pl]} {\capit S. Plewik,} {\it On completely Ramsey sets,}
Fundamenta Mathematicae, vol. 127 (1986), pp. 127-132.
\smallskip
\itemitem{[Sa]} {\capit G. Sacks,} {\it Forcing with perfect closed sets,}
in: Axiomatic Set Theory, Proc. Symp. Pure Math., vol. 13 (1971),
pp. 331-355.
\smallskip
\itemitem{[Sh]} {\capit S. Shelah,} {\it Proper and improper
forcing,} in preparation.
\smallskip
\itemitem{[Ve]} {\capit B. Veli{\v c}kovi\'c,} {\it CCC posets
of perfect trees,} Compositio Mathematica, vol. 79 (1991),
pp. 279-294.
\smallskip